\magnification=1200

\def\la#1{\hbox to #1pc{\leftarrowfill}}
\def\ra#1{\hbox to #1pc{\rightarrowfill}}
\def\fract#1#2{\raise4pt\hbox{$ #1 \atop #2 $}}
\def\decdnar#1{\phantom{\hbox{$\scriptstyle{#1}$}}
\left\downarrow\vbox{\vskip15pt\hbox{$\scriptstyle{#1}$}}\right.}

\parskip=1pc
\def\za{\vrule height6pt width4pt depth1pt}

\def\ds{\displaystyle}
\def\lrar{\longrightarrow}
\def\ext{\Lambda}
\def\power{\Gamma}

\def\tensor{\otimes}
\def\deg{\hbox{deg}}
\def\dim{\hbox{dim}}
\def\sp#1{SP^{#1}}
\def\spy{SP^{\infty}}
\def\hol#1{\hbox{Hol}_{#1}}

\def\map#1{\hbox{Map}_{#1}}

\font\twelverm=cmr10 at 12pt
\font\ninerm=cmr9
\font\eightrm=cmr8
\font\sc=cmcsc10
\font\svntnrm=cmr17


\centerline{\svntnrm  The Geometry of the Space of Holomorphic}
\medskip
\centerline{\svntnrm
Maps from a Riemann Surface to a}
\medskip
\centerline{\svntnrm  Complex Projective Space}

\vskip20pt

\centerline{\twelverm Sadok
Kallel\footnote{\raise3pt\hbox{$\dag$}}{\eightrm The author 
holds a Postdoctoral fellowship with CRM, Montr\'eal}}

\centerline{\ninerm\it Centre de Recherches Math\'ematiques,
Universit\'e de Montr\'eal} 
\vskip 10pt
\centerline{\twelverm R. James
Milgram\footnote{\raise3pt\hbox{$\dag\dag$}}{\eightrm Research
partially supported by a grant from the NSF and a 
grant from the CNRS}}

\centerline{\ninerm\it Stanford University, California.}
\vskip 20pt

\noindent{\bf\twelverm\S0. Introduction}

In recent years there have been a number of papers on the homology
and geometry of spaces of holomorphic maps of the Riemann sphere into complex
varieties, [S], [C$^2$M$^2$], [MM1], [MM2], [BHMM], [Gu].  However,
the very classic question of the structure of the spaces of holomorphic
maps from complex curves $M_g$ of genus $g\geq 1$, to complex varieties has
proved to be very difficult.  There is Segal's stability theorem, [S],
which shows that 
the natural inclusion of the space of based holomorphic maps of degree $k$ into
the space of all based maps, $Hol_k^*(M_g, V) \hookrightarrow Map_k^*(M_g, V)$
is a homotopy equivalence through a range of dimensions which increases with
$k$ when $V = {\bf P}^n$, the complex projective space.  Also, there is the
extension of this result by J. Hurtubise to further $V$; [H]. But that is
about all.

In this paper we begin the detailed study of the topology of
the $Hol_k^*(M_g, {\bf P}^n)$.  We are able to completely
determine these spaces and
their homology when $M_g$ is an elliptic curve and we give an
essentially complete determination in the case where $M_g$ is
hyperelliptic.  In particular we determine the rational homology
of these spaces when $k \ge 2g-1$ in the elliptic and hyperelliptic
cases.

Let $M_g$ be a genus $g$ complex curve.  The key analytic result on
the structure of $Hol_k^*(M_g, {\bf P}^1)$ is
Abel's theorem which identifies the disjoint pairs of $k$-tuples
of unordered points $\langle r_1, \dots, r_k\rangle$ and $\langle
p_1, \dots, p_k\rangle$ on $M_g$ that are the roots
and poles of a meromorphic function on $M_g$ in terms of the
Abel-Jacobi embedding  of
$M_g$ into its Jacobi variety, $\mu\colon M_g \ra{1} J(M_g)$.  This
extends to give
necessary and sufficient conditions for an $(n+1)$ tuple $\{V_i~|~V_i
 = \langle
x_{i,1}, \dots ,x_{i,k}\rangle, i = 0, \dots, n\}$ of points in $M_g$
with $\cap_0^n \{V_i\} = \emptyset$ to be the root data for
a holomorphic map of $M_g$ into ${\bf P}^n$, $n\ge 2$, and thus
defines an embedding $Hol_k^*(M_g, {\bf P}^n) \subset (SP^k(M_g))^{n+1}$,
where $SP^k(X)$ is the $k$-fold symmetric product of $X$.

We begin by studying a compactification
of the space $Hol_k^*(M_g,{\bf P}^n)$ which we denote $E_k^n(M_g)$
obtained
by taking the closure of the embedding
$$Hol_k^*(M_g,{\bf P}^n)
\subset (SP^k(M_g))^{n+1}$$
above.  We show in (2.2)
that $E_k^n(M_g)$
is given as the total space of a fibering,
$$ ({\bf P}^{k-g})^{n+1} \ra{1.5} E_k^n(M_g) \ra{1.5} J(M_g)$$
for $k\ge 2g-1$, with $H^*(E_k^n(M_g);{\bf A}) = H^*({\bf P}^{k-g},{\bf A})^{n+1}\otimes
H^*(J(M_g);{\bf A})$ for ${\bf A}$ any commutative ring of coefficients.  
Moreover, for $k \le 2g-2$, $E_k^n(M_g)$ is stratified by strata 
which are fiberings
$$({\bf P}^{s})^{n+1} \ra{1} S_k^s \ra{1} A_k^s$$
where the $A_k^s$
are subspaces of $J(M_g)$ determined by the curve.  Specifically,
$W_k^j\subset J(M_g)$ is the subspace of points in the image of the
$k^{th}$ Abel-Jacobi map
$$\mu_k\colon SP^k(M_g) \ra{1} J(M_g)$$
with
$(\mu_k)^{-1}(x) = {\bf P}^s$ where $s \ge \hbox{max}(k-g,0)+j$, 
and $A_k^j = W_k^j - W_k^{j+1}$.   
Both spaces $W_k^j$ and $A_k^j$ are extensively
studied in [ACGH], [G] (cf. \S5).

This compactification has the property that $Hol_k^*(M_g, {\bf P}^n)$ is
open in $E_k^n(M_g)$ and we write $V_k^n(M_g)$ for the (closed)
complement $E_k^n(M_g)-Hol_k^*(M_g, {\bf P}^n)$.
The space $E_k^n(M_g)$ being a closed, compact manifold, 
Alexander-Poincar\'e duality now gives
$${\tilde H}^{2k(n+1) - 2ng - *}(Hol_k^*(M_g, {\bf P}^n);{\bf F})
\cong H_*(E_k^n(M_g), V_k^n(M_g);{\bf F})$$
for $k \ge 2g-1$.
This then indicates that one can understand $Hol_k^*(M_g, {\bf P}^n)$
by first studying $V_k^n(M_g)$.

The space $V_k^n(M_g)\subset E_k^n(M_g)$ is the union of two pieces.
The first, for $k > 2g-1$, is a 
subfibration $Z_k^n(M_g)$ of the form
$$\left\{\bigcup_{i=0}^{n+1}({\bf P}^{k-g})^i\times {\bf P}^{k-g-1}
\times({\bf P}^{k-g})^{n-i}\right\}
\ra{1.5} Z_k^n(M_g) \ra{1.5} J(M_g)$$
with analogous definitions of $Z_k^n(M_g)$ for $k \le 2g-1$ (see
Definition 2.6).  Moreover,
it is easy to determine the relative cohomology groups,
$H^*(E_k^n(M_g), Z_k^n(M_g);{\bf F})$ for $k >2g-1$,
and not too difficult for $k \le 2g-2$, provided we know
enough about the $W_k^j$.

In the case of hyperelliptic and elliptic curves the
Riemann-Roch theorem determines
the $W_k^j$ explicitly.  
Moreover, a complete
description of the $W_k^j$ for all curves of genus $\le 6$ is
given in [ACGH, p. 206 - 211].

In \S5~we study the way
in which the $W_k^i$ determine $H^*(E_k^n(M_g), Z_k^n(M_g);{\bf A})$
in detail and obtain a spectral sequence converging to these
groups with explicit $E^1$-term.

\noindent{\bf Proposition} 5.9:  \tensl Suppose that $k \le g$, then
there is a spectral sequence converging to $H^*(E_k^n(M_g),
Z_k^n(M_g);{\bf A})$ 
with $E^1$-term
$$E^1 = H_*(SP^k(M_g), SP^{k-1}(M_g); {\bf A})
\oplus \coprod_iH_*(W_k^i, W_{k-1}^i;{\bf A})
\otimes
{\tilde H}^*(S^{2i(n+1)};{\bf A})$$
and in case $k \ge g$ then
$$E^1 = \Sigma^{2(k-g)(n+1)}H_*(J(M_g), W_{k-1}^{k-g};{\bf A}) \oplus
\coprod_{i> k-g}H_*(W_k^i, W_{k-1}^i;{\bf A})
\otimes {\tilde H}_*(S^{2i(n+1)};{\bf A}).$$
\tenrm

The second part of $V_k^n(M_g)$ consists of the image of
an action map
$$\nu\colon~ M_g\times E_{k-1}^n(M_g) \ra{1} E_k^n(M_g)$$
introduced in (2.4) that puts in redundant roots,
and is the main source of difficulty in recovering
the homology of $Hol_k^*(M_g, {\bf P}^n)$ from the homology groups
$H_*(E_k^n(M_g); {\bf A})$ or $H_*(E_k^n(M_g), Z_k^n(M_g);{\bf A})$.
In order to handle this part we construct a spectral
sequence which converges to the
homology of the pair $(E_k^n(M_g), V_k^n(M_g))$ starting with
the relative groups
$$H_*(E^n_k(M_g), Z^n_k(M_g);{\bf F}).$$

\noindent{\bf Theorem} 4.3:  \tensl There is a spectral sequence
converging to $H_*(E_k^n(M_g), V_k^n(M_g);{\bf F})$ for any field ${\bf F}$
with $E^1$-term
$$\eqalign{E^1~=~\coprod_{i+j=k\atop i\ge 1} H_*&(E_i^n(M_g), Z_i^n(M_g);{\bf F})\otimes
H_*(SP^j(\Sigma M_g), SP^{j-1}(\Sigma M_g); {\bf F})\cr
&\oplus H_*((SP^k(\Sigma M_g), SP^{k-1}(\Sigma M_g); {\bf F})\cr}$$
\tenrm

The rest of \S4~ discusses the structure of $d_1$ and various
properties such as the multiplicative pairing of the spectral
sequences discussed immediately after the proof of (4.3).  Of
course the structure of the groups $H_*(SP^j(\Sigma M_g),
SP^{j-1}(\Sigma M_g); {\bf F})$
is well known from e.g. [DT], [D], and [M1].  It is also reviewed
in \S6.

Next we use these results to clarify the structure of the natural
inclusions
$$Hol_k^*(M_g, {\bf P}^n) \hookrightarrow Map_k^*(M_g, {\bf P}^n).$$ 
The spectral sequence of (4.3) has a natural break at $i = 2g-1$
in the sense that $i \ge 2g-1$ implies that the relative homology
groups
$$H_*(E_i^n(M_g), Z_i^n(M_g); {\bf F}) = \cases{0& $* < 2(i-g)(n+1)$,\cr
H_{*-2(i-g)(n+1)}(J(M_g); {\bf F})& otherwise\cr}$$
and, when $i< 2g-1$, the
groups depend on the structure of $\mu_i$, the $W_i^r$, and have
to be determined case by case.  We call the case $i \ge 2g-1$
the stable range for the spectral sequence and completely determine
the differentials in this range.  For $* > (2n-1)(k-2g+1)$ only
the stable homology contributes to $H_*(E_k^n(M_g), V_k^n(M_g);{\bf F})$.  
Thus, by Poincar\'e duality, the unstable range only contributes
homology above this range.  On the other hand it is easily seen that
the duals of the stable range classes inject into
$H_*(Map_k^*(M_g, {\bf P}^n);{\bf F})$ and many of these classes live
considerably above the stable range above.  Thus, as was the case 
with $Hol_k^*({\bf P}^1, {\bf P}^n)$ ([C$^2$M$^2$]) one has considerably 
more information about the map than was given by Segal's stability theorem.

In the final sections we study the cases of elliptic and
hyperelliptic curves.  Here, as indicated,
the Riemann-Roch theorem gives complete
control of the $W_k^r$ and consequently the $E^1$-term of the
spectral sequence is within range of calculation.  In particular,
for elliptic curves the situation is completely understood.

Typical of the results in the elliptic case is

\noindent{\bf Lemma} 9.4:  \tensl Let $I = (b_1, b_2)$ be the
augmentation ideal in the polynomial ring ${\bf Q}[b_1, b_2]$ where
$dim(b_i) = 2$, $i = 1, 2$.  Then
$$H_*(Hol_k^*(M_1, {\bf P}^1);{\bf Q}) \simeq \left\{{\bf Q}[b_1, b_2]/I^{k+1}
\right\}(1, 
e_1, e_2, h_1, h_2, v) \oplus {\bf Q}(w_1, w_2, \dots, w_{2k-1})$$
where $dim(e_i) = 1$, $dim(h_i) = 2$, $dim(w_i) = 2k-3$ and
$dim(v) = 3$.
\tenrm

Here $M_1$ is an arbitrary elliptic curve.  Additionally, it
turns out that the map
$$H^*(Map_k^*(M_1, {\bf P}^1); {\bf Q}) \ra{1.5} H^*(Hol_k^*(M_1, {\bf P}^1);
{\bf Q})$$
is surjective for all $k \ge 1$ in this case.

For hyperelliptic curves there are some technical questions which
seem difficult to handle for finite field coefficients, but with some
effort a complete determination of the $E^1$-term in (4.3) with
${\bf Q}$-coefficients
is given in \S10 --15, together with sufficient differentials
to completely determine $E^{\infty}$ and conclude that $H_*(Hol_k^*(M_g,
{\bf P}^n); {\bf Q})$ injects into $H_*(Map_k^*(M_g, {\bf P}^n);{\bf Q})$ for
$k \ge 2g-1$ in the hyperelliptic case as well.

As the arguments are pretty involved we summarize the salient points
here.

To begin, the Riemann-Roch theorem gives a complete determination of the
$W_i^r$ as quotients of $SP^i(M_g)$ via an action discussed in the
proof of (10.2)
$$SP^r({\bf P}^1) \times SP^l(M_g) \ra{1.5}  SP^{l+2r}(M_g)$$
induced from the
Abel-Jacobi map $\mu_2\colon SP^2(M_g) \ra{1} J(M_g)$ which fails
to be an embedding at only one point where it has inverse image
a copy of ${\bf P}^1$,
the ${\bf P}^1$ in the action map above.

These observations give the following result for
hyperelliptic curves:

\noindent{\bf Lemma}  10.2:  \tensl  Suppose that $M_g$ is 
hyperelliptic and $\tau \in J(M_g)$ is the hyperelliptic point.
Let $k\geq 1$ and $t\leq \left[ {k\over 2}\right]$. Then
\item{(a)}  for $k \le g$, the space $W_k^t$
is $\mu_{k-2t}(SP^{k-2t}(M_g)) + t\tau$,
\item{(b)}  for $2g-1 > k > g$, $t>k-g$, we also have
$W_k^t=\mu_{k-2t}(SP^{k-2t}(M_g)) + t\tau$.
\tenrm

Of course, this gives the $W_k^t$ as quotients, so, in order to
obtain information about the $W_k^t$ we introduce some spectral
sequences which take care of the details of the quotienting process
in \S11~through \S13.  Using them we are able to determine the
rational homology of the $W_k^t$ as follows.

\noindent{\bf Lemma} 13.7:  \tensl 
\item{(a)}  The inclusion $W_j \subset J(M_g)$ induces an injection in
rational homology $H_*(W_j;{\bf Q}) \hookrightarrow H_*(J(M_g);{\bf Q}) =
\Gamma(e_1, \dots, e_{2g})$ with image the subvector space spanned by the
subspaces
$$\left\{ \Gamma_s(e_1, \dots, e_{2g})[M_g]^t~|~s+t \le j\right\}$$
where $[M_g]= \sum_1^ge_{2i-1}e_{2i}$ is the image of the fundamental
class of $M_g$ under the Abel-Jacobi map $\mu_*$.
\item{(b)}  $H_*(W_{j-1}; {\bf Q})$ injects into $H_*(W_j;{\bf Q})$ under the
inclusion so the relative groups are given as
$$H_*(W_j,W_{j-1};{\bf Q}) \cong H_*(W_j;{\bf Q})/H_*(W_{j-1};{\bf Q}).$$
\tenrm

Next we turn to the homology of the spaces $Hol_k^*(M_g, {\bf P}^n)$
in the hyperelliptic cases.  This involves using the spectral sequences
and calculations above.  Of course, the $E^1$-term in the spectral sequence
(4.3) is very complex,
but one is able to recognize in it the direct sum of a family
of chain complexes, each of which calculates a part of a certain $Tor$
or $Ext$ group of the exterior algebra $\Gamma(e_1, \dots, e_{2g})$
modulo the two sided ideal generated by $[M_g]$.  This explains
\S14~ which is devoted to the calculation of the relevant $Tor$-groups.

Finally, in \S15~ we are able to put these results together to obtain
our main calculational result,

\noindent{\bf Theorem} 15.1:  \tensl The natural map $H_*(Hol_k^*(
M_g, {\bf P}^n); {\bf Q}) \ra{1} H_*(Map_k^*(M_g, {\bf P}^n);{\bf Q})$ is
injective for $k\ge 2g$ and $n>2$ if $M_g$ is hyperelliptic.\tenrm

\noindent{\sc Acknowledgements:}
We would like to thank the ``Centre de Recherches Math\'ematiques"
at the University of Montr\'eal for their support and hospitality 
during the period when much of this work was done. 
The special case of (based) meromorphic functions on the torus; i.e.
$g=1, n=1$ was worked out as part of the first author's Ph.D thesis
written under the direction of the second author.
We also thank Professor 
J. Hurtubise for numerous comments and his generous support.


\noindent{\bf\twelverm\S1. Preliminaries On the Abel-Jacobi Map}

We review some classical definitions and theorems about the algebraic
geometry of curves.  We begin by defining the Abel-Jacobi map together
with the Jacobi variety $J(M_g)$ associated to any positive genus Riemann
surface.  Good references are [ACGH] and [G].

Any Riemann surface $M_g$ of genus $g\ge 1$ has
$g$ independent holomorphic sections of the cotangent bundle $\tau^*(M_g)$,
(holomorphic 1-forms), $w_1, w_2,\dots, w_g$, the {\tensl abelian
differentials} on $M_g$.

Fix a basepoint $p_0\in M_g$.  Then we can associate to each point $p\in M_g$
and each path $\gamma$ between $p_0$ and $p$ the vector of integrals
$$\mu_{\gamma}(p) = \left(\int_{p_0}^pw_1, \dots, \int_{p_0}^pw_g\right)\in
{\bf C}^g.$$
Since any two paths $\gamma$ and $\gamma^{\prime}$ between $p_0$ and $p$ together
determine a closed loop based at $p_0$, (which we will denote as $L$),
we have that $\mu_{\gamma}(p)$
is well defined up to vectors of the form
$$\left(\int_Lw_1, \dots, \int_Lw_g\right)$$
If we choose a set of loops $L_1,\dots ,L_{2g}$
which, in homology, form a basis for $H_1(M_g;{\bf Z})$ they give rise to the
following $g\times 2g$ {\tensl period matrix}
$$\Omega = \pmatrix{\int_{L_1}w_1& \cdots & \int_{L_{2g}}w_1\cr
\vdots &\ddots& \vdots\cr
\int_{L_1}w_g&\cdots & \int_{L_{2g}}w_g\cr}$$
Thus, since the $w_i$ are {\tensl closed}, the values $\mu_{\gamma}(p)$ depend
only on $p$ and not $\gamma$ in the quotient torus
$$J(M_g)~=~{\bf C}^g/\Omega \cong (S^1)^{2g}.$$
Consequently, they give rise to a well defined map $\mu \colon M_g 
\ra{1} J(M_g)$ which is called the {\tensl Abel-Jacobi map} for
$M_g$.

At this point we need to introduce symmetric products.

\noindent{\bf Definition}:  The {\tensl $n$-fold symmetric product of the
space $X$},$\sp{n}(X)$, is defined to
be the set of unordered $n$-tuples of points of $X$; i.e.
$\sp{n}(X) = X^n/{\cal S}_n$ where ${\cal S}_n$ is the symmetric group
on $n$ letters. A point in $\sp{n}(X)$ will be written in the form
$\sum m_i x_i$ with $m_i>0$ and $\sum m_i = n$, or in the form
$$\langle x_1, x_2, \dots, x_n\rangle.$$

\noindent{\bf Remarks}:  
Let $X$ be any $CW$ complex with base point $*$, then there are inclusions
$\sp{n} (X)\hookrightarrow \sp{n+1} (X)$
which identify $\sum_in_iP_i$ with
$\sum_in_iP_i+*$, and we get the increasing sequence
of spaces
$$* = \sp{0} (X)\subset\sp{1} (X)\subset\cdots\subset\sp{n-1} (X)\subset
\sp{n} (X)\subset\cdots.$$
The union of this sequence is the infinite symmetric product $\spy(X)$,
based at $*$.  
Moreover, if $X$ is path connected, then the homotopy type of $\spy(X)$
is independent of the choice of $*$.

We can  extend the Abel-Jacobi map to the symmetric products 
of $M_g$ by the rule $\langle m_1,\dots,m_k
\rangle \mapsto \mu(m_1)+ \cdots + \mu(m_k)$, obtaining the family
of maps
$$\mu_k \colon SP^k(M_g) \ra{1.5} J(M_g).\leqno{1.1}$$
The map $\mu$ is called the Abel-Jacobi map and it, together with the
extensions $\mu_k$ are a critical piece of the structure data for $M_g$.

It is a remarkable result due to Andreotti that for any $M_g$ the
symmetric products $SP^k(M_g)$ are all complex manifolds, indeed,
complex algebraic varieties of complex dimension $k$. To see this
note that the symmetric product $\sp{k}({\bf C})$ is diffeomorphic to
${\bf C}^k$ via the map that takes
the unordered collection $\langle z_1, \dots, z_k\rangle$ of points in
${\bf C}$ to 
the coefficients of the monic polynomial of degree $k$ with the $z_i$,
$1\le i \le k$ 
as roots.  Since locally $\sp{k}(M_g)$ is modeled on $\sp{k}({\bf
C})$, the result follows. 

\noindent{\bf Remark} 1.2:  In the special case that $M_g ={\bf P}^1$ a
slight extension 
of the above argument identifies $\sp{k}({\bf P}^1)$ with ${\bf P}^k$
which is now
regarded as the space of all homogeneous polynomials of degree $k$ in
$2$ variables. 

\noindent{\bf Remark} 1.3:  In (1.14) we point out an extension of this
result due to Mattuck  which, for $k\ge 2g-1$ identifies $\sp{k}(M_g)$ with
the total space of a fibration over $(S^1)^{2g}$ with fiber ${\bf P}^{k - g}$.

\noindent{\bf Remark} 1.4: Another way of thinking about $J(M_g)$
is as the Picard group of $M_g$, the
space of isomorphism classes of holomorphic line bundles on $M_g$.  From this
point of view the addition in $J(M_g)$ corresponds to the tensor product
of line bundles.  Under this correspondence, the map $\mu$ takes $m\in M_g$
to the line bundle obtained from the trivial bundle over $M_g$ by gluing in
a copy of the negative Hopf bundle over $m$ (see [G]).

We can associate to every holomorphic function $f\in Hol(M_g, {\bf P}^1)$
a divisor 
$(f)$ defined by $(f) = \sum n_i Z_i - \sum m_j P_j$ where $\{Z_i\}$,
$\{P_i\}$ are 
the zeros and poles respectively of $f$ with multiplicities $n_i$ and
$m_i$ respectively.   
By standard residue calculations, it is easy to see that $\sum n_i=\sum m_j$.
The space of pairs of disjoint divisors on
a Riemann surface 
$(\zeta, \eta)$, subject to the condition $\hbox{deg}(\zeta) =
\hbox{deg}(\eta)$ 
constitutes then the first step in the description
of the space of holomorphic maps from
the surface to the Riemann sphere ${\bf P}^1$.

For maps of ${\bf P}^1$ to itself this description is enough
for there do exist meromorphic maps having prescribed roots and
poles (of equidegree); i.e. any divisor $D=\zeta-\eta$, 
$\deg D=0$ is the divisor of a meromorphic function on ${\bf P}^1$.
For the general case $g\ge 1$ it turns out that a pair $(\zeta, \eta)$
as above need not necessarily give rise to a meromorphic function and
one needs a further condition.

\noindent{\sc Theorem} 1.5~(Abel):  \tensl Given positive divisors $D$ and
$D^{\prime}$ on $M_g$, $\deg D=\deg D^{\prime}$, then there exists an $f\in
Hol(M_g,{\bf P}^1)$ 
so that $(f) = D- D^{\prime}$ if and only if $\mu(D) = \mu(D^{\prime})$.
\tenrm

The $f$ associated to the difference $D-D^{\prime}$ is unique provided
we specify in advance the image of $p_0$ (based maps) and neither
$D$ nor $D^{\prime}$ contains the basepoint, $p_0$.  Additionally,
given $D-D^{\prime}$, there are unique {\sl disjoint} positive divisors
$D_1$, $D_1^{\prime}$ so that $D-D^{\prime} = D_1 - D_1^{\prime}$ and
any $f$ with $(f) = D-D^{\prime}$ will have roots precisely the terms in
$D_1$ and poles in $D_1^{\prime}$.  We make the following definition. 

\noindent{\bf Definition} 1.6:
The {\tensl divisor} space $\hbox{Div}_k(X)$ for a given space
$X$ is the set of pairs of {\tensl disjoint} positive divisors on $X$, i.e.
$$\hbox{Div}_k(X) = \left\{(D,D^{\prime})\in SP^k(X)\times SP^k(X)~|~
D\cap D^{\prime} = \emptyset\right\}
\leqno{1.7}$$
and more generally
$$\eqalign{\hbox{Div}^n_k(X) = \left\{ (D_1, \dots, D_{n+1}~|~ 
D_i \in SP^k(X),\right. &~1\le i \le n+1,\cr
&\left. D_1\cap D_2 \cap \cdots \cap D_{n+1} = \emptyset\right.\}.\cr}
\leqno{1.8}$$

\noindent{\sc Corollary} 1.9:  \tensl The space of based holomorphic maps of
degree $k$, $Hol^*_k(M_g,{\bf P}^1)$, is the
inverse image of $~0$ under the subtraction map
$$s\colon \hbox{Div}_k(M_g-p_0) \ra{1.5} J(M_g)$$
where the subtraction map $s$ is given by $s((D, D^{\prime})) = \mu(D)-
\mu(D^{\prime})$.\tenrm

More generally, and perhaps more usefully, we have

\noindent{\sc Corollary} 1.10:  \tensl The space of based holomorphic maps of
degree $k$, $Hol^*_k(M_g, {\bf P}^n)$ is the subspace of
$\hbox{Div}^n_k(M_g-p_0)$ consisting of $(n+1)$-tuples of
degree $k$ positive divisors, subject to the following constraint:
$$\mu_k(D_1) = \mu_k(D_2) = \cdots = \mu_k(D_{n+1}).$$ 
\tenrm

It was using this formulation of the space $Hol^*_k(M_g, {\bf P}^n)$ that
Segal [S] proceeded to prove his stability result. We will shortly give and use
a yet more explicit version of this corollary. But before we do so,
here are some of the standard results on the Abel-Jacobi map that we will be using
\item{(1.11)} The orginal map $\mu\colon M_g \ra{1} J(M_g)$ is an
embedding and the (complex) dimension of the image of $\mu_d\colon SP^d(M_g)
\ra{1} J(M_g)$ is $d$ for $d\le g$.  In particular $\mu$ is onto for 
$d\ge g$ {\tensl (Jacobi inversion theorem)}.
\item{(1.12)} The preimage of any point $\mu (p)\in J(M_g)$,
$\mu_d^{-1}(\mu (p))\in SP^d(M_g)$ is always a complex
projective plane ${\bf P}^m$ for some $m\ge 0$.
\item{(1.13)} For $d\le 2g-2$ the dimension $m$ is less than or equal to
${d\over 2}$ {\tensl (Clifford)}.
\item{(1.14)} For $d\ge 2g-1$ the map $\mu_d$ makes $SP^d(M_g)$ into an
analytic fiber bundle over $J(M_g)$ with fiber ${\bf P}^{d-g}$ {\tensl(
Mattuck [Mt])}.

In \S9 we will give more details on the structure of these maps for
$g \le 5$, but now we turn to the construction of an explicit model
for the space $\hol{k}^*(M_g,{\bf P}^n)$ from our considerations thus far.


\noindent{\bf\twelverm\S2. A Compactified Version of
$\hol{k}^*(M_g,{\bf P}^n)$} 

It should be apparent from \S1~that we are interested in
the inverse image of $0$ under $s\colon\hbox{Div}_k(M_g-x_0)\rightarrow 
J(M_g)$ for holomorphic maps to ${\bf P}^1$,
and generally in the inverse image of the iterated diagonal
$$\Delta^{n+1}(J(M_g) \subset \underbrace{J(M_g)\times \cdots\times J(M_g)}_{
n+1-times}$$
under $\mu_k\times \cdots \times \mu_k$ for maps into ${\bf P}^n$ .

Let the space $E_{i_0,i_1,\ldots, i_n}$ be
defined as the fiber product of $\mu\times\cdots\times \mu$ and $\Delta$
in the diagram below
$$\matrix{E_{i_0,i_1,\ldots, i_n}&\fract{}{\ra 2}&\sp{i_0}(M_g)
\times\cdots\times\sp{i_n}(M_g)\cr
\decdnar{}&&\decdnar{\mu\times\cdots\times \mu}\cr
J&\fract{\Delta}{\ra 2}&J(M_g)\times\cdots\times J(M_g).\cr}
\leqno{2.1}$$
More explicitly
$$E_{i_0,i_1,\ldots, i_n} = \{(D_0, \ldots, D_n)\in
\sp{i_0}(M_g)\times\cdots\times\sp{i_n}(M_g)~|~
\mu (D_i)=\mu (D_j)\}.$$
Clearly $\hol{k}^*(M_g, {\bf P}^n)\subset E_{k,k,\dots, k}$ as the 
subset consisting of all the $\{ (D_0,\dots, D_n)\}$ where no $D_i$ contains 
$*$ and $\bigcap_0^n D_i = \emptyset$.  Of course, for $k$ sufficiently 
small $\hol{k}^*(M_g,{\bf P}^n)$ is empty and this is consistent
with the fact that Riemann surfaces of positive genus do
not admit holomorphic functions with a single pole.  
However, once $\hol{k}^*(M_g, {\bf P}^n)$ is non-empty,
and $k$ is sufficiently large, $\hol{k}^*(M_g, {\bf P}^n)$ will be open
in $E_{k,\dots, k}$ with $E_{k,\dots, k}$ as its closure.  
Then $E_{k,k,\dots, k}$
is the compactification of $\hol{k}^*(M_g, {\bf P}^n)$ mentioned in the title 
of this section.

\noindent{\bf Lemma} 2.2:  {\sl~
For $i_j\geq 2g-1$, $0\le j \le n$, we have a fibering
$$\ds\prod_{j=0}^n{\bf P}^{i_j-g}\lrar E_{i_0,i_1,\ldots, i_n}\lrar J(M_g)$$
and the fiber is totally non-homologous to zero so
$$H_*(E_{i_0,i_1,\ldots, i_n})\cong H_*({\bf P}^{i_0-g})\tensor \cdots\tensor
H_*({\bf P}^{i_n-g})\tensor 
H_*(J(M_g)).$$}
\noindent{\sc Proof:} In this range of dimensions $i_j\geq 2g-1$, one 
can see by virtue of Mattuck's theorem that the space $E_{i_0,i_1,\ldots, i_n}$
becomes the total space of a fibration
$${\bf P}^{i_0-g}\times\cdots\times{\bf P}^{i_n-g}\lrar E_{i_0,i_1,\ldots, i_n}
\lrar J(M_g)$$
obtained from the pull-back of a product of Mattuck's fibrations
$${\bf P}^{i_0-g}\times\cdots\times{\bf P}^{i_n-g}\lrar 
\sp{i_0}(M_g)\times\ldots\times\sp{i_n}(M_g)\lrar J(M_g)^{n+1}$$
by the diagonal inclusion $\Delta$. Since $\Delta$ is injective and since the 
Serre spectral sequence for Mattuck's
fibration collapses at $E_2$, the lemma follows from standard spectral
sequence comparison arguments.
\hfill\za

\noindent{\tensl Further Properties of the $E_{i_0, \dots, i_n}$}

First of all, we notice that for any $j\geq 0$, we have natural inclusions 
$$E_{i_0,i_1,\ldots, i_n}
\hookrightarrow E_{i_0,\ldots,i_j+1,\ldots, i_n}$$
given by adding the basepoint in the $(j+1)^{st}$ position.

Secondly we observe that the map
$\mu\times \cdots \times \mu\colon
SP^{i_0}(M_g)\times \cdots \times SP^{i_n}(M_g) \ra{1} J(M_g)^{n+1}$
is multiplicative by construction and so it induces a multiplicative
pairing on the $E_{i_0, \dots, i_n}$ which is commutative and associative:
$$\nu^{I,J}\colon E_{i_0,\dots, i_n}\times E_{j_0,\dots, j_n} \ra{1.5}
E_{i_0+j_0, \dots, i_n + j_n}.
\leqno{2.3}$$
Also, $E_{1,1,\dots, 1} = M_g$ and the pairing above thus yields an action
$$\nu\colon \left(\coprod_{k=1}^{\infty} \sp{k}(M_g)\right)\times 
E_{i_0, \dots, i_n}
\ra{1.5} \coprod_{k=1}^{\infty}E_{i_0+k, \dots,i_n+k}
\leqno{2.4}$$
explicitly defined via the diagonal multiplication
$$\sp{k}(M_g)\times\left( SP^{i_0}(M_g)\times\cdots\times SP^{i_n}(M_g)\right)
\ra{1.5} SP^{i_0+k}(M_g)\times \cdots \times SP^{i_n+k}(M_g),$$
which acts on points as follows
$$\eqalign{(\sum m_s, (\langle m_{11},\dots, m_{1,i_0}\rangle,
& \dots , \langle m_{n1},\dots,m_{ni_n}
\rangle))\cr
& \mapsto (\langle \sum m_s, m_{11},\dots,m_{1,i_0}\rangle,\dots,
\langle \sum m_s, m_{n1},\dots,m_{ni_n}\rangle.\cr}$$
We can then give an explicit reformulation of the description of 
$\hol{i}^*(M_g, {\bf P}^n)$ in these terms.

\noindent{\bf Lemma} 2.5:  \tensl Let ${\tilde E}_{i_0,\dots,
i_n}\subset E_{i_0,\dots, 
i_n}$ be the subspace in which no $m_{rs}$ is equal to the base point
$*$.  Then, 
$$ Hol^*_i(M_g, {\bf P}^n)\cong {\tilde E}_{i,i,\dots, i} - \hbox{Image}(\nu)
\cap{\tilde E}_{i,i,\dots, i}.$$
\tenrm

The following constructions are now needed for the remainder of our discussion.

\noindent{\bf Definition} 2.6:  \tensl For each $n$, the space $LE_i$
is the quotient 
$$LE_i~=~ E_{i,i,\dots, i}\big\slash \left\{\bigcup E_{i,i,\dots, i-1,
i,\dots, i}\right\}$$ 
The space $QE_{i}$ is the quotient
$$LE_i\big\slash \left\{\hbox{Image}(\nu)\right\}.$$
\tenrm

\noindent{\bf Remark} 2.7: When $i\geq 2g-1$, the space
$E_{i,\ldots, i}$ is a manifold 
of dimension $2(n+1)i-2ng = 2(i-g)(n+1)+2g$. However when $i<2g-1$
there is no garantee that $E_{i,\ldots, i}$ is actually a manifold.

Using Alexander-Poincar\'e duality, we can now deduce that

\noindent{\bf Lemma} 2.8:\tensl~Assume $i\geq 2g-1$. Then for
untwisted coefficients  
${\bf A}$ we have
$$H_j(Hol_i^*(M_g, {\bf P}^n);{\bf A})\cong H^{2(i-g)(n+1)+2g -
j}(QE_i;{\bf A}).$$
\tenrm

Note at this point that the multiplication $\mu$ of (2.3) induces
an associative,  
commutative multiplication on the $QE_i$ as well:
$$\mu\colon QE_i\times QE_j \ra{1.5} QE_{i+j}.$$
From lemma (2.8), it is clear that it is the space $QE_i$
that we wish to study
in the remainder of this paper.  Unfortunately, it is generally very hard to 
obtain the cohomology of such a space without a careful analysis of
the piece we 
collapse out.  So, in order to do this we follow the
procedure of [C$^2$M$^2$] and replace the cone on the union above by a
more complex but much more structured space.


\noindent{\bf\twelverm\S3.  A Model for $QE_i$}

Consider the twisted product space
$$DE(M_g) = 
\left(\bigcup 
E_{i_0,i_1,\ldots, i_n}\right)\times_t\spy (cM_g),
\leqno{3.1}$$
where $cT$ denotes the reduced cone on $T$ and where the twisting
$t$ is given by the action above.  Precisely, points of $DE$ are
of the form 
$$\{(D_0,\ldots, D_n),(t_1,z_1)\ldots (t_l,z_l)\},~D_i\in\sp{k_i}(M_g)~
\hbox{and}~\mu (D_i)=\mu (D_j)\}$$
with the identification that when $t_i=0$ the point above is identified with
$$\{(D_0+z_i,\ldots, D_n+z_i),(t_1,z_1)\ldots {\widehat{(t_i,z_i)}}\ldots
(t_l,z_l)\}$$
where the entry $(t_i,z_i)$ is deleted from the last set of coordinates.
Clearly $\mu (D_r+z_i)=\mu (D_s+z_i)$ and the construction makes sense.

The space $DE$ is filtered by the subspaces
$$DE_{k_0,k_1,\ldots, k_n}(M_g) =
\bigcup_{{{i_0+l\leq k_0}\atop\vdots}\atop {i_n+l\leq k_n}}
E_{i_0,i_1,\ldots, i_n}\times_t\sp{l} (cM_g).$$
Observe that there are projection maps
$$p_{k_0,k_1,\ldots, k_n}\colon
DE_{k_0,k_1,\ldots, k_n}\lrar E_{k_0,k_1,\ldots, k_n}/
\left\{\hbox{Image}(\nu)\right\}$$
where
$$(v_1, \dots, v_s, (t_1, w_1), \dots, (t_r, w_r))\mapsto \left\{(v_1, \dots,
v_s)\right\}$$
and the inverse images of points consist of
contractible sets. This implies that
the maps $p_{k_0,k_1,\ldots, k_n}$
are acyclic and induce isomorphisms
$$H_*(DE_{k_0,k_1,\ldots, k_n};{\bf F} )\lrar 
H_*(E_{k_0,k_1,\ldots, k_n}/\left\{\hbox{Image}(\nu)\right\};{\bf F}).$$
We combine this with the isomorphism in lemma 2.8 to get

\noindent{\bf Corollary} 3.2: {\sl~
Let $k\ge 2g-1$, then
$$\ds \tilde H^{2k(n+1)-2ng-*}(\hol{k}^*(M_g,{\bf P}^n); {\bf F} )\cong
H_*(DE_{\underbrace{k,\ldots,k}_{n+1}}/
\bigcup_i DE_{k,\ldots,\underbrace{k-1}_{ith-entry},\ldots,k}; {\bf F} )$$}

This result is very useful because it is possible, using the
filtration of the space $DE_{k,\dots, 
k}$ described above, to construct spectral
sequences with 
known $E_1$-terms which converge to the cohomology of the relative spaces above
for all $k\ge 1$.


\noindent{\bf\twelverm\S4.  The Spectral Sequence}

\noindent{\tensl The diagonal action:}
The diagonal multiplication introduced in 2.3
$$\sp{r}(M_g)\times E_{i_0,i_2,\ldots,i_n}\fract{\nu}{\ra 3}
E_{i_0+r,\ldots,i_n+r}$$
induces an action in homology, $\nu_*$,
$$\nu_*\colon \left(\coprod_{r=0}^{\infty}H_*(\sp{r}(M_g);{\bf F})\right)
\otimes \left(\coprod_{j=1}^{\infty}H_*(LE_i;{\bf F})\right)\ra{1.5}
\coprod_{j=1}^{\infty}H_*(LE_{i+r}; {\bf F})\leqno{4.1}$$
for any field coefficients ${\bf F}$.
Of course, this quotient action fits together with the original action on the
$E_{i,\dots, i}$ via the following commutative diagram
$$\matrix{H_*(\sp{r}(M_g);{\bf F})\otimes H_*(E_{i,\dots, i};{\bf F})
&\fract{id\otimes p_*}{\ra{1.5}}&
H_*(\sp{r}(M_g);{\bf F})\otimes H_*(LE_i;{\bf F})\cr
\decdnar{\nu_*}&&\decdnar{\nu_*}\cr
H_*(E_{i+r, \dots, i+r};{\bf F})&\fract{p_*}{\ra{1.5}}&
H_*(LE_{i+r};{\bf F}).\cr}\leqno{4.2}$$

\noindent{\bf Theorem} 4.3: \tensl
There is a spectral sequence converging to $H_{*}(QE_{k}; {\bf F})$
with $E^1$ term
$$\eqalign{E^1~& =~ \coprod_{i+j = k\atop i \ge 1}
{\tilde H}_*(LE_i ; {\bf F})
\tensor H_*(\sp{j}(\Sigma M_g),\sp{j-1}(\Sigma M_g);{\bf F})\cr
&\oplus
H_*(\sp{k}(\Sigma M_g), \sp{k-1}(\Sigma M_g);{\bf F})\cr}$$
and $d_1( \Theta \otimes \{ |a_1| \cdots |a_r|\}) = (-1)^{||\Theta|| ||a_1||}
\nu_*(a_1\otimes \Theta)\otimes \{|a_2| \cdots |a_r|\}$.
\tenrm
 
\noindent{\sc Proof:}
From the model for $QE$ constructed in \S3 we have that
$$QE_k\cong
DE_{k,\ldots,k}/
\bigcup DE_{k,\ldots,k-1,\ldots,k}$$
where
$$DE_{k,\ldots, k}(M_g) =
\bigcup_{i_r+j\leq k}E_{i_0,i_1,\ldots, i_n}\times_t\sp{j} (cM_g).$$
Also, $LE_i=E_{i,\dots, i}\big\slash \left\{\bigcup E_{i,i,\dots, i-1,
i,\dots,i}\right\}$ 
so we can write
$$QE_k = \bigcup_{i + j = k}LE_i\times_t\left(
SP^j(cM_g)/SP^{j-1}(cM_g)\right)$$ 
where the twisting in the description above is given as before by
$$(h, \{(0,w_1), (t_2, w_2), \dots, (t_r, w_r)\}) \sim (\nu(w_1, h),
\{(t_2, w_2), \dots, 
(t_r, w_r)\}).$$
To obtain the desired spectral sequence, introduce the filtration by
$${\cal F}_r(QE_k) ~=~ \bigcup _{i\le r\atop i+j=k}LE_i\times_t\left(
SP^j(cM_g)/SP^{j-1}(cM_g)\right),$$
and the remainder of the proof of the theorem is direct.\hfill\za

\noindent{\tensl A multiplicative structure for the spectral sequence}:  The
induced multiplication on the $QE_i$'s,
$${\bar \nu}\colon QE_i\times QE_j \ra{1.5} QE_{i+j}$$
passes to the spectral sequences above and defines a pairing of $E^1$-terms:
$$\eqalign{&\left( H_*(LE_i;{\bf F})\otimes H_*(SP^j(\Sigma M_g),
SP^{j-1}(\Sigma M_g);{\bf F})\right)\cr 
&\hbox to .3in{\hfill} \otimes
\left( H_*(LE_v;{\bf F})\otimes H_*(\sp{w}(\Sigma M_g),
\sp{w-1}(\Sigma M_g);{\bf F})\right)\cr 
&\hbox to .7in{\hfill}
\ra{1.5} H_*(LE_{i+v};{\bf F})\otimes H_*(\sp{w+j}(\Sigma M_g),
\sp{w+j-1}(\Sigma M_g);{\bf F})\cr}$$ 
for which the $d_i$'s act as derivations.  For this reason it is often
convenient to consider 
all the spectral sequences above at once.  In particular we can
describe the direct sum 
of all the $E_1$-terms as a trigraded ring where an element $x\in E_1$
has tridegree $(i,j,*)$ if and only if
$$x\in H_s(LE_i)\tensor H_{*-s}(\sp{j}(\Sigma M),\sp{j-1}(\Sigma
M);{\bf F} ).$$ 

\noindent{\bf Remark} 4.4:  
There is a related spectral sequence for $H_*(QE_k;{\bf F})$ obtained
by filtering
$$QE_k = \bigcup_{i + j = k}LE_i\times_t\left(
SP^j(cM_g)/SP^{j-1}(cM_g)\right)$$ 
in a somewhat different way.  Instead of filtering by $i$ in the
expression above, filter by the number of {\tensl distinct} $t$'s in the point
$(h, \{(t_1,w_1), (t_2, w_2), \dots, (t_r, w_r)\})$.
In the space $SP^{\infty}(\Sigma M_g)$ this filtration results in the
Eilenberg-Moore 
spectral sequence with $E_2$-term $Ext_{H_*(\spy(M_g);{\bf
F})}^{*,*}({\bf F}, {\bf F})$. 
To describe the resulting spectral sequence most efficiently it is
best to include all the 
$QE_k$'s at once, and what results in our case is a trigraded $E_2$-term,
$$Ext_{H_{*,*}(\spy(M_g);{\bf F})}^{r,s,*}({\bf F},
\bigoplus_{i=1}^{\infty}H_*(LE_i;{\bf F}))$$ 
where the summand such that $r+s=k$ corresponds to 
the $E_2$-term of the spectral sequence converging to $H^*(QE_k;{\bf F})$.

\noindent{\tensl Some remarks on differentials.}
There are similar spectral sequences for the $Div$-spaces (cf. 1.6)
starting (as in 3.1) with the model
$$\left(\coprod_{j=0}^{\infty}\underbrace{\sp{j}(M_g)\times 
\cdots\times \sp{j}(M_g)}_{n+1-times}\right)\times_T\spy(cM_g)
\leqno{4.5}$$
where $T$ is now the diagonal twisting which identifies the point 
$(0,z)$ in the cone $cM$ with the diagonal element $\Delta^{n+1}(z)$
in $(M_g)^{n+1}$ and then extends this multiplicatively.
The associated spectral sequences have $E_1$-term
$$\coprod_{i+j=k}\left[H^*(\sp{i}(M_g), \sp{i-1}(M_g);{\bf F})\right]^{n+1}
\otimes H^*(\sp{j}(\Sigma M_g), \sp{j-1}(\Sigma M_g);{\bf F})\leqno{4.6}$$
and converge to $H^*(Div_{k,\dots, k}(M_g);{\bf F})$.  (Here the
superscript $n+1$ 
means the $(n+1)$-fold tensor product.)  Moreover, as is clear, the
inclusions
$$E_{i,\dots, i} \subset \sp{i}(M_g)\times \cdots \times \sp{i}(M_g)$$
induce maps of spectral sequences here (in cohomology)
to the spectral sequences above for the $\hol{k}^*$ spaces.  
(Or in homology from the spectral sequences for the $\hol{k}^*$ spaces
to these for the $Div$-spaces.) 

But for the spectral sequences for the $Div$-spaces, and using H. Cartan's
little constructions [Car] to embed the homology into the chain complex one is
able to construct 
an explicit (small) filtered chain complex with associated spectral
sequence equal to that in $(4.3)$, [K1].
In particular one knows that $E_{\infty} = E_1$ for $n\ge 2$ for the
$Div$-spaces spectral sequence 4.6.

In the next section we will identify a region of the spectral sequence for the $\hol{k}^*$
spaces where the induced map of homology spectral sequences is an injection.  Hence,
in this range for $n>1$ the spectral sequence for $\hol{k}^*(M_g, {\bf P}^n)$
collapses at $E_1$.  When $n=1$ there are differentials, however our knowledge of
the differentials in the $Div$-spectral sequence here implies considerable information
about the differentials for the $\hol{k}^*$ spaces here as well.

\noindent{\tensl The $d^1$-differential for the highest filtration terms of the spectral
sequence.}  One region where the two sequences don't compare very well is the tail end
of (4.3), the terms
$${\tilde H}_*(M_g;{\bf F})\otimes H_*(\sp{k-1}(\Sigma M_g), \sp{k-2}(
\Sigma M_g);{\bf F})~\hbox{and }~H_*(\sp{k}(\Sigma M_g),
\sp{k-1}(\Sigma M_g);{\bf F}).$$ 
More generally it can happen that $\mu_k\colon SP^k(M_g) \ra{1} J(M_g)$ is an
embedding for $1\le k \le m(M_g)$ in which case we have

\noindent{\bf Lemma} 4.7:  \tensl If the Abel-Jacobi map
$\mu_k\colon SP^k(M_g) \ra{1} J(M_g)$ is an embedding for $1\le k\le m(M_g)$,
then in this range we have $E_{k,\dots, k} \cong SP^k(M_g)$,
$E_{k,\dots, k-1, k,\dots,k} 
\cong SP^{k-1}(M_g)$ included in $E_{k,\dots, k}$ via the usual
base-point embedding 
$$SP^{k-1}(X)\subset SP^k(X),\hbox to .3in{\hfill}
\langle x_1, \dots, x_{k-1}\rangle \mapsto
\langle x_1,\dots, x_{k-1},*\rangle$$
and the action $\nu$ corresponds to the usual multiplication
$$SP^r(M_g)\times SP^{k-r}(M_g) \ra{1.5} SP^k(M_g).$$
Consequently, in this range we have $LE_{k,\dots, k} \cong
SP^k(M_g)/SP^{k-1}(M_g)$, 
$1\le k \le m(M_g)$, and the spectral sequence in this region is the
corresponding spectral 
sequence for the quasi-fibration $SP^{\infty}(M_g) \ra{1}
SP^{\infty}(cM_g) \ra{1} SP^{\infty}(\Sigma M_g)$.
\tenrm

\noindent(The only statement above which might need clarification is
the last.  But 
recall that the spectral sequence with field coefficients has $E_2$-term
$$\bigoplus_{k,t} H^*(SP^k(M_g), SP^{k-1}(M_g);{\bf F})
\otimes H^*(SP^t(\Sigma M_g), SP^{t-1}(\Sigma M_g);{\bf F})$$
and all the differentials preserve the sum $k + t$.  Moreover, in this
region the 
chain embedding techniques of [K1] are valid, so the two spectral
sequences have the same (internal) differentials.)

For example, if $M_g$ is not hyperelliptic then $\mu_2\colon SP^2(M_g)\ra{1}
J(M_g)$ is an embedding.  Also the map $\mu_3\colon SP^3(M_5)\ra{1} J(M_5)$
is an embedding for most curves of genus 5.  (See the discussion in
[ACGH], chapter V.)


\noindent{\bf\twelverm\S5.  The Jacobi varieties $W_i^j$ and
a spectral sequence for the $LE_i$ spaces}

Our next step is to analyze the spectral sequence of 4.3.
To do this we must understand
the groups $H_*(LE_i)$. The $LE_i$ are quotients of the 
spaces $E_i\equiv E_{i,\ldots,i}$ defined in \S2~and these latter spaces turn
out to be built out of fibrations with projective spaces as fibers.  Here
there is a stratification of the image of $SP^n(M_g)$ in $J(M_g)$
and over each stratum we get such a fibration, though the dimensions
of the fibers vary as we move from stratum to stratum.

\noindent{\sc Definition} 5.1:  \tensl The image of $\mu_d$
in $J(M_g)$ is written $W_d$. Also the set of points $y\in W_d$
so that $\mu_d^{-1}(y) = {\bf P}^m$ with $m \ge r$
is denoted $W_d^r$ so we have a decreasing filtration
$$W_d \supseteq W_d^1 \supseteq W_d^2 \supseteq \cdots \supseteq W_d^r\supseteq
\cdots\,\supseteq W_d^{[d/2]+1}=\emptyset .$$
\tenrm

It is well known that $W_d=J(M_g)$ whenever $d\geq g$, and that the
dimension $m$ of a {\sl generic} fiber ${\bf P}^m$ over $W_d$
is $0$ when $d\leq g$ and $d-g$ when $d\geq g$.

{\parskip=2pt
\noindent{\sc Examples}:
\item{(5.2)} The map $\mu_1$ is always an embedding and so $W_1\cong M_g$.
\item{(5.3)} In the genus 1 case the original map $\mu = \mu_1$ is
the identity
and $J(M_1) = M_1$ while the $\mu_d$ are fiberings for $d\ge 2$.}

\noindent{\tensl Assume now that $g\ge 2$.}
\item{(5.4)} The map $\mu_2\colon SP^2(M_g) \ra{1} J(M_g)$ is an
embedding unless $M_g$ is hyperelliptic (cf \S10).  In the hyperelliptic case
$W_2^1$ is a single point $p$ and $\mu^{-1}_2(p) = {\bf P}^1$.  It follows
that $W_2$ can be identified with $SP^2(M_g)$ with a single ${\bf P}^1$
blown down.  Indeed, we can check that the normal bundle to ${\bf P}^1\subset
SP^2(M_g)$ is $\xi^{1-g}$, the line bundle with self-interesection number
$1-g$, and from this it follows that
\itemitem{(1)} $SP^2(M_2)$ is the blow up of $J(M_2)$ at a single
point (here any genus 2 surface is automatically hyperelliptic). 
\itemitem{(2)} For $g\ge 3$ we have that $W_2 = SP^2(M_g)/({\bf P}^1\sim *)$ is
a manifold with a single isolated point singularity which looks like the
cone on a Lens space $L_{g-1}^3$.
\item{}

We can introduce the complementary 
subspace $A_k^i = W_k^i-W_k^{i+1}$. By definition,
we have that $\forall x\in A_k^i$, $\mu^{-1}_k(x)={\bf P}^i$. The result
of Mattuck quoted in (1.14) takes actually the more general form

\noindent{\bf Theorem} 5.5:~
(Mattuck)~\tensl
$\mu^{-1}_k(A_k^i)\subset\sp{k}(M_g)$ is the total space 
of a locally trivial analytic fibration
${\bf P}^i\rightarrow\mu^{-1}_k(A_k^i)\rightarrow A_k^i.$
\tenrm

\noindent{\bf Remark} 5.6:  The following formula relating $W_k^i$ and
$W_{k-1}^{i-1}$ is a special case of [G, (20), p. 53]:
$$W_k^i~=~\bigcap_{m \in \mu(M_g)}(W_{k-1}^{i-1} + m),\leqno{(**)}$$
that is, $W_k^i$ is obtained as the intersection of all translates of
elements of $W_{k-1}^{i-1}$ by elements of $\mu(M_g)\subset J(M_g)$.
In particular, this shows that 
$$W_k^i \subset W_{k-1}^{i-1}+* \subset W_k^{i-1}.\leqno{(***)}$$
From now on we do not differentiate between $W_{k-1}^{i-1}\subset W_{k-1}$
and $W_{k-1}^{i-1}+*\subset W_k$; denoting both by $W_{k-1}^{i-1}$.

Furthermore, Gunning, [G, p. 54], gives the following dimension estimates
for the associated containments.

\noindent{\bf Lemma} 5.7:  \tensl If the subvariety $W_{k-1}^{i-1}$ is
non-empty then 
\item{$\bullet$}  dim $W_k^i < $dim $W_{k-1}^{i-1}$ whenever $2 \le i
\le k \le g$. 
\item{$\bullet$}  dim $W_{k-2}^{i-1} < $ dim $W^{i-1}_{k-1}$ whenever
$2 \le i \le k \le g$.
\tenrm

\noindent{\bf Remark} 5.8:  In chapter V of [ACGH] techniques for
determining the $W_k^i$'s are
extensively discussed and results for $g\le 6$ are completely given,
(p. 206-211).

We now turn to the pull-back spaces $E_{k,\ldots,k}$ and observe
that the fibration of lemma (5.5) induces in turn a fibration
$$({\bf P}^i\times\cdots\times{\bf P}^i)\lrar {\cal A}_k^i
\fract{\pi}{\ra 2} A_k^i,~{\cal A}_k^i\subset E_{k,\ldots,k}.$$
When collapsing ${\cal A}_{k-1}^i$, we get a quotient
$X={\cal A}_k^i/{\cal A}_{k-1}^{i}\hookrightarrow E_{k,\ldots,k}/
E_{k-1,\ldots,k-1}$ and a quotient map $X\lrar LE_{k,\ldots,k}$.
The space $X$ projects down via $\mu$ to
$$A_k^i/A_{k-1}^i= W_k^i/W_{k-1}^i;$$ 
this last equality being a consequence of the inclusion
$W_{k}^{i+1}\subset W_{k-1}^i$ described above.

We then pass to the quotient $LE_k$ and observe that since
$W_k^i\subset W_{k-1}^{i-1},~i\geq 1$, we must further collape,
along fibers this time, subsets of the form
$$\bigcup_{j\leq n+1}({\bf P}^i)^{j-1}\times{\bf P}^{i-1}\times({\bf
P}^i)^{n+1-j} \hookrightarrow ({\bf P}^i)^{n+1}.$$
The fiber $({\bf P}^i)^{n+1}$ has a top $2i(n+1)$ cell of the form
$e^{2i(n+1)}=
e^{2i}_1\times\cdots\times e_{n+1}^{2i}$ where $e^{2i}_j$ is the top $2i$
dimensional cell in the $j$-th copy ${\bf P}^i$ with boundary mapping to
${\bf P}^{i-1}\subset{\bf P}^i$. We then see that
$$({\bf P}^i)^{n+1}/\cup_j({\bf P}^i)^{j-1}\times{\bf
P}^{i-1}\times({\bf P}^i)^{n+1-j}
\simeq e^{2i(n+1)}/\partial e^{2i(n+1)}\simeq S^{2i(n+1)}$$
The space so obtained is denoted by $T_k^i$. In the case when $i=0$
we have that $W_k^1\subset W_{k-1}\subset W_k$, and hence
$W_k/W_{k-1}\cong\sp{k}(M_g)/\sp{k-1}(M_g)$.

Now we filter the Jacobian according to the increasing sequence of $W_i$'s
$$J(M_g)=W_g\supset W_{g-1}\supset\cdots\supset W_1=M_g.$$ 
This induces a filtration on $LE_k$ yielding a spectral sequence which
by the preceeding discussion has $E^1$ term as follows.

\noindent{\bf Proposition} 5.9:~\tensl Suppose that $k\leq g$, then
there is a spectral sequence converging to ${\tilde H}_*(LE_{k})$ with
$E^1$ term 
$$E^1= H_*(\sp{k}(M_g),\sp{k-1}(M_g))\oplus
\coprod_i H_*(W_k^i,W_{k-1}^i)\tensor{\tilde H}_*(S^{2i(n+1)})$$
and in case $k\geq g$, then
$$E^1= \Sigma^{2(k-g)(n+1)}H_*(J(M_g),W_{k-1}^{k-g})
\oplus
\coprod_{i>k-g}H_*(W_k^i,W_{k-1}^i)\tensor{\tilde H}_*(S^{2i(n+1)}).$$
\tenrm

\noindent{\bf Remark} 5.10: Note that $H_*(LE_k)= H_*(J(M_g))
\tensor{\tilde H}_*(S^{2(k-g)(n+1)})$ whenever $k> 2g-1$ for then
$W_{k-1}^{k-g}=\emptyset$ and $W_k^i=W_{k-1}^i=J(M_g)$ otherwise.

\noindent{\bf Remark} 5.11: The spectral sequence of (5.9) turns
out to collapse at $E^1$ for all cases we treat in this paper.


\noindent{\bf\twelverm\S6. The Structure of Symmetric Products}

In this section, we describe the homology of the symmetric products
$\sp{n}(M_g)$ for all $n$ and $g$.  Also, since it is required in the
spectral sequences of (4.3) and (5.9), we give the structure
of $H_*(SP^n(\Sigma M_g))$
again for all $n$, $g$.  In the case of $SP^n(M_g)$ we follow the
description given by I.G. Macdonald, [Mc], while the description
for the suspension is taken from [M1] and unpublished work of N.
Steenrod.

There is an evident pairing
$\sp{n} (X)\times \sp{m} (X)\fract{+}{\ra 1}\sp{n+m} (X)$ given by
addition of points and
this turns $\spy (X)$ into an associative, abelian
monoid with $*$ as a two sided identity.  
The Dold and Thom states that
$$\spy(X) \simeq \prod_1^{\infty} K(H_i(X;{\bf Z}), i)
\leqno{6.1}$$
is a product of Eilenberg-MacLane spaces if $X$ is path connected
[DT]. 
Applying this to $X=M_g$ yields
$$\spy (M_g)\simeq K({\bf Z}^{2g},1)\times K({\bf Z},2)
\simeq (S^1)^{2g}\times {\bf P}^{\infty} \leqno{6.2}$$
where $K({\bf Z},1) \simeq S^1$ while
$K({\bf Z},2) \simeq {\bf P}^{\infty}$, the infinite complex projective space.
We then find that
$$H_*(\spy (M_g);{\bf Z})\cong\ext (e_1,\dots, e_{2g})\tensor
\power (a)\leqno{6.3}$$
where $\ext(,)$ denotes the exterior algebra on the stated generators
while $\power(a)$ denotes the {\tensl divided power algebra} on $a$:
the ring with ${\bf Z}$-generators $a_i$ $i = 0,1, \dots$ and multiplication
$a_i a_j = {i+j\choose j}a_{i+j}$.
This accords well with the multiplication of the $\sp{n}(M_g)$ described
above. Often however it is convenient to work with the cohomology
rings so we need the following description.

\noindent{\tensl I.G. Macdonald's description of $H^*(\sp{k}(M_g);{\bf Z})$}

Consider the map $[M_g]^* \colon M_g \ra{1} {\bf P}^{\infty} = K({\bf Z}, 2)$,
taking the fundamental class to the dual of the orientation class, and the
map $\vee e_i\colon M_g \ra{1} K({\bf Z}^{2g}, 1) \simeq (S^1)^{2g}$.
Both ${\bf P}^{\infty}$ 
and $(S^1)^{2g}$ are associative abelian $H$-spaces with the structure
on ${\bf P}^{\infty}$ 
coming from the identification $\spy({\bf P}^1) = {\bf P}^{\infty}$
described earlier. 
  
This allows us to extend $[M_g]^*$ and $\vee e_i$ to a multiplicative
map
$$\theta(k)\colon (M_g)^k \ra{1} {\bf P}^{\infty}\times (S^1)^{2g}$$
which, by definition factors through the symmetric product $SP^k(M_g)$.  
(In the limit, as $k\mapsto \infty$ this gives the Dold-Thom equivalence 
$\spy(M_g) \ra{1} K(({\bf Z})^{2g},1)\times K({\bf Z}, 2)$.)

With respect to the maps $\theta(k)$ we can describe the cohomology
ring $H^*(\sp{k}(M_g); {\bf Z})$ as follows:

\noindent{\bf Theorem} (I.G. Macdonald) 6.4:  \tensl The cohomology
ring of $\sp{k}(M_g)$ 
over the integers ${\bf Z}$ is generated by the elements
$$f_1= e_1^*, \dots , f_i= e_{2i-1}^*, \dots f_g=e_{2g-1}^*,
f_1^{\prime} = e_2^* , \dots f_{g}^{\prime}=
e_{2g}^*,~\hbox{and}~b$$
subject to the following relations:
\item{(i)}  The $f_i$'s and the $f_i^{\prime}$'s anti-commute with
each other and 
commute with $b$;
\item{(ii)} If $i_1, \dots ,i_a, j_1,\dots ,j_b,k_1,\dots, k_c$ are
distinct integers 
from 1 to $g$ inclusive, then
$$ f_{i_1}\cdots f_{i_a}f_{j_1}^{\prime} \cdots
f_{j_b}^{\prime}(f_{k_1}f_{k_1}^{\prime} 
- b) \cdots (f_{k_c}f_{k_c}^{\prime} - b) b^q = 0$$
provided that 
$$a+b+2c+q = n+1.$$

If $k < 2g$ all the relations above are consequences of those for
which $q=0$, and 
if $n > 2g-2$ all the relations are consequences of the single relation
$$b^{k-2g+1}\prod_{i=1}^g(f_if_i^{\prime}-b) = 0.$$
\tenrm

\noindent(Actually Macdonald only stated this result in [Mc] for fields of
characteristic zero as coefficients, however, since the relevant invariant maps
that he used to prove the result are actually surjective over ${\bf Z}$ it is
quite direct to extend the result to integers as well.)

\noindent{\tensl The homology of $\sp{n}(X)$ for more general $X$}

The homology of the spaces $\sp{n}(X)$ splits according to
a result of Steenrod, [D], (see also [M1]) which holds for $X$ any
$CW$-complex: 
$$H_*(\spy (X);{\bf Z}) = \bigoplus_jH_*(\sp{j}(X),\sp{j-1}(X);{\bf Z} ).
\leqno{6.5}$$
Consequently, the ring $H_*(\spy (M_g);{\bf Z})$
is bigraded, by defining $x\in H_*(\spy(X); {\bf Z})$ to have bidegree
$(i,j)$ iff 
$x\in H_j(\sp{i}(M_g),\sp{i-1}(M_g);{\bf Z} )\subset H_*(\spy
(M_g);{\bf Z})$.   
The bigrading is multiplicative in the sense that the product
map
$$\sp{m}(X)\times \sp{n}(X) \ra{1} \sp{n+m}(X)$$
in homology induces
a bigraded ring map
$$H_{i,j}(\sp{n}(X);{\bf A})\otimes H_{r,s}(\sp{m}(X); {\bf A}) \ra{1.5}
H_{i+r, j+s}(\sp{n+m}(X);{\bf A}),$$
where the first degree is the dimension and the second the bidegree.

The bigrading is preserved by the diagonal
map so the cohomology ring is bigraded as well, and the two structures
together form a bigraded Hopf algebra.  In the case of the $\sp{k}(M_g)$ we
have that the coproduct on each of the $1$-dimensional generators is primitive
while the $2$-dimensional generator $[M_g]$ in $H_{1,2}(\spy (M_g);{\bf Z})$
corresponding to the orientation class of $M_g$ has coproduct
$$\Delta ([M_g]) = [M_g]\otimes 1 + \sum_1^g (e_{2i-1}\otimes e_{2i} -
e_{2i}\otimes e_{2i-1}) + 1\otimes [M_g].\leqno{6.6}$$

\noindent{\tensl The homology of $\sp{n}(\Sigma (M_g))$}

We will need
the cohomology and homology of the spaces $\sp{n}(\Sigma (M_g))$ where 
$\Sigma X$ denotes the suspension of $X$.  From the Dold-Thom theorem
$$\spy(\Sigma M_g) \simeq \left(\prod_1^{2g}{\bf P}^{\infty}\right)\times
\spy(S^3)\leqno{6.7}$$ 
where we can assume that the generator for the homology of each ${\bf
P}^{\infty}$
has bidegree $(1,2)$.  Thus,
this amounts to describing the
homology and cohomology of $K({\bf Z},3)$ and its associated bigrading.

With coefficients ${\bf F}_p$ for any odd prime $p$  we have, [M1],
$$\eqalign{
H_*(K({\bf Z},3);{\bf F}_p)~&\cong \ext (|a|, 
|\gamma_p|, \cdots \, |\gamma_{p^i}| \cdots)
\cr
& ~\otimes \Gamma (|a^{p-1}|a|) \otimes \Gamma (|\gamma_p^{p-1}|\gamma_p|) \otimes
\cdots \otimes \Gamma (|\gamma_{p^i}^{p-1}|\gamma_{p^i}|) \otimes \cdots\cr}
\leqno{6.8}$$
The generator $\gamma_{p^i}$ corresponds to
$[\sp{p^i} (M_g)]$ and $|\gamma_{p^i}|$
has filtration degree $p^i$ and dimension $2p^i+1$ while
$|\gamma_{p^{i-1}}^{p-1}|\gamma_{p^{i-1}}|$ has bidegree
$(p^i,2p^i + 2)$.

In the case of the prime $2$ we have
$$H_*(K({\bf Z},3); {\bf F}_2)~\cong~\power[|a|, |a_2|, \dots, |a_{2^i}|, ...].
\leqno{6.9}$$
\noindent{\bf Remark} 6.10:  With ${\bf F}_p$ coefficients,
the divided power algebra $\power [a]$ splits as
an algebra into a tensor product of truncated polynomial algebras [Car]
$$\power[a] \cong {\bf F}_p[a]/a^p\otimes {\bf
F}_p[a_p]/(a_p)^p\otimes \cdots$$ 
 
When ${\bf F}={\bf Q}$, divided power algebras are isomorphic to
polynomial algebras and the answer simplifies greatly; namely
$$H_*(\spy (\Sigma M_g);{\bf Q} )= {\bf Q} [|e_1|,\ldots, |e_{2g}|]\tensor
\ext (a).\leqno{6.11}$$

The answer is simpler to describe if we use cohomology.  Here,
$H^*(\spy(S^3);{\bf F}_p)$ 
is always a tensor product of an exterior algebra on odd dimensional generators
and a polynomial algebra on the even dimensional generators.  In other words
in the formulae 6.8 above one replaces all the divided power algebras by
polynomial algebras with generators of the same bidegrees to get the cohomology
ring description.

We take advantage of the bigrading to write $H_*(\spy (\Sigma
M_g);{\bf F}_p)$ in the form 
$$\coprod_1^{\infty}
H_*(\sp{n}(\Sigma M_g), \sp{n-1}(\Sigma M_g); {\bf F}_p).$$

These summands can be written out more completely in terms of the bigrading of
$\spy (S^3)$ as follows.
First, we write $H_*(\spy(S^3);{\bf F}_p)=H_*(\spy(\Sigma S^2);{\bf F}_p)$
in this way.  
The decomposition here has an alternate description, [C$^2$M$^2$]:
$$H_*(\spy(S^3); {\bf F}_p) = \bigoplus_1^{\infty} \Sigma^{4i}{\cal D}_i^*(p)
\leqno{6.12}$$
where the ${\cal D}_i$ are the Snaith splitting components of the
loop space $\Omega^2 S^3$ for the prime $p$ and ${\cal D}_i^*$ means
dual, where we index the dual by $\dim({\cal D}_{i,j}) = -j$.

Then, going back to (6.7) we see that we can write the term
$$H_*(\sp{n}(\Sigma M_g), \sp{n-1}(\Sigma M_g);{\bf F}_p) =
\bigoplus_{j=0}^n\Sigma^{4j}
{\cal D}_j^*(p)\otimes {\bf F}_p[b_1, \dots, b_{2g}]_{n-j}
\leqno{6.13}$$
and ${\bf F}_p[b_1,\dots, b_{2g}]_{n-j}$ is the free ${\bf Z}$-module on the
degree $n-j$ monomials in the variables $b_1, \dots, b_{2g}$.  There
are ${n-j + 2g-1\choose 2g-1}$ such monomials.

\noindent{\bf Example}:  What follows is a list of generators
for $H_{*}(\sp{j}(\Sigma M),\sp{j-1}(\Sigma M);{\bf Z}_p )$
in the $E^1$ term mod($p$) of the spectral sequence of (4.3)
together with their tridegrees 
$$ \matrix {Generator && trigrading\cr
|M|&&(0,1,3)\cr
|e_i|&&(0,1,2) \cr
|\gamma_{p^i}|&&(0,p^i, 2p^i+1)\cr
|\gamma_{p^{i-1}}^{p-1}|\gamma_{p^{i-1}}|
&&(0,p^i, 2p^i+2).\cr}
\leqno{6.14}$$

\noindent{\tensl Explicit $d^1$-differentials in the spectral sequence.}

In 4.3 the differential $d^1$ is implicitly determined.  Now that we
have the explicit form of the homology groups $H_*(SP^n(M_g))$ and
$H_*(SP^n(\Sigma M_g))$ we can make $d^1$ explicit.  For example we have
$$d^1(|f_1|^{s_1}\cdots |f_{2g}|^{k-s_1 - \cdots -s_{2g-1}})
 = \sum_1^{2g}f_i\otimes|f_1|^{s_1}\cdots |f_i|^{s_i-1}\cdots |f_{2g}|^{k-s_1
- \cdots - s_{2g-1}}
\leqno{6.15}$$
and
$$\eqalign{d^1(|M_g||f_1|^{s_1}\cdots|f_{2g}|^{k-s_1-\cdots -1}~& =~
[M_g]\otimes 
|f_1|^{s_1}\cdots|f_{2g}|^{k-s_1 - \cdots -1}\cr
&- \left(d^1(|f_1|^{s_1}\cdots |f_{2g}|^{k-s_1-\cdots -1}\right)
\left(1\otimes |M_g|\right).\cr}
\leqno{6.16}$$
It is easily checked that this part of the $d^1$-differential is
injective to the term 
${\tilde H}^*(M_g;{\bf F})\otimes {\bf F}[|f_1|, \dots, |f_{2g}|]_{k-1}$
with cokernel 
spanned by the monomials of the following form:
$$\{ f_t|f_i|^{j_i}|f_{i+1}|^{j_{i+1}}\cdots |f_{2g}|^{k-1-j_1 -\cdots
- j_{2g-1}}\},~~\hbox{where}~ j_i> 0,~\hbox{and}~t> i. 
\leqno{6.17}$$

\noindent{\bf Remark} 6.18:  This is the only time the $d^1$ differential
on the elements $|f_1|^{i_1}\cdots |f_{2g}|^{i_{2g}}$ is non-trivial,
since in the remaining parts
of the spectral sequence their images lie in the part which has been
collapsed out. 

\noindent{\bf Example} 6.19:  
If $g=1$ then the differentials in this region have the
form $d^1(|f_1|^s|f_2|^{k-s}) = f_1\otimes|f_1|^{s-1}|f_2|^{k-s}
+ f_2\otimes|f_1|^s|f_2|^{k-s-1}$
provided $s$ and $k-s$ are both greater than zero and
$$\eqalign{d^1(|T||f_1|^s|f_2|^{k-s-1}) ~&=~ f_1f_2\otimes
|f_1|^s|f_2|^{k-s-1}\cr & - f_1\otimes |T||f_1|^{s-1}|f_2|^{k-s-1}
 - f_2\otimes |T||f_1|^s|f_2|^{k-s-2}.\cr}$$
The role of this last differential can be regarded as identifying terms of the
form $f_1f_2\otimes |f_1|^s|f_2|^{k-s-1}$ with terms involving $|T|^*$.  Also,
the first of these differentials is injective but not surjective.  The
quotient has as a basis the set of images of the elements
$f_2\otimes |f_1|^s|f_2|^{k-s-1}$ with $s> 0$.  Consequently it has dimension
$k-2$.


\noindent{\bf\twelverm\S7.  The Stable Range}

In the case where $i \ge 2g$, we've seen already that the Abel-Jacobi
map becomes an analytic fibration, and hence
the projection $E_{i,\ldots, i}\ra{1} J(M_g)$ fibers as
$${\bf P}^{i-g}\times\ldots\times{\bf P}^{i-g}\lrar E_{i,\ldots,
i}\lrar J(M_g). \leqno{7.1}$$
Because of this and other stabilization properties (\S 8), we refer to the 
range $i \ge 2g$ as {\tensl stable}.
In this case, we see that the relative groups
$${\tilde H}(LE_i, {\bf A}) \cong H_*(E_{i,i,\ldots,i}, 
\bigcup E_{i, \dots, i-1, \dots, i}; {\bf A})$$
are given by
$$\eqalign{\left[H_*({\bf P}^{i-g}, {\bf P}^{i-g-1};{\bf A})\right]^{n+1}&
\otimes H_*(J(M_g);{\bf A})\cr
& \cong \Sigma^{2(i-g)(n+1)}H_*(J(M_g);{\bf A}).\cr}
\leqno{7.2}$$
Dually, in cohomology the relative group is given as
$$b_0^{i-g}b_1^{i-g}\cdots b_n^{i-g}\otimes H^*(J(M_g);{\bf A})=
(b_0\cdots b_n)^{i-g}\otimes H^*(J(M_g);{\bf A}),
\leqno{7.3}$$
a form which is often of more use in calculations.

\noindent{\bf Remark} 7.4:
Of course, in both these forms, the inclusions ${\bf P}^{i-g}\subset
\sp{i}(M_g)$ 
induce injections in homology onto ${\bf Z}$-direct summands, but even more,
from Macdonald's results (6.4) it follows that the inclusion of pairs
$$({\bf P}^{i-g}, {\bf P}^{i-g-1})\subset (\sp{i}(M_g),\sp{i-1}(M_g))$$
induces an inclusion in homology sending $H_{2(i-g)}({\bf P}^{i-g},
{\bf P}^{i-g-1};{\bf Z})$ 
isomorphically to the group $H_{2(i-g)}(\sp{i}(M_g),
\sp{i-1}(M_g);{\bf Z}) = {\bf Z}$. 
(We will expore this fact further below.)

It follows then that in this stable
range the $E_1$-term for $\hol{k}^*(M_g, {\bf P}^n)$ injects into the
corresponding 
$E_1$-term for the $Div$-space as was asserted at the end of \S4.
In particular 
this gives as a corollary to the results of [K1],

\noindent{\bf Lemma} 7.5:  \tensl The spectral sequence of (4.3)
collapses at $E_1$ in the stable range for $n\ge 2$.
\tenrm

Now we turn to some fairly direct calculations which lead to the determination
of the $d^1$-differential in the $\hol{k}^*$ spectral sequence for $n=1$ in the
stable range.

More exactly,  in the stable range the relative group injects,
$$H^*(E_{i,i},  E_{i,i-1}\cup E_{i-1,i}; {\bf A})
\subset H^*(E_{i,i}, {\bf A})
$$
as the principal ideal
$$((b_0b_1)^{i-2g}W_0W_1)
\subset H^*(E_{i,i}; {\bf A})
\leqno{7.6}$$
where $W_j$ is the polynomial
$$\eqalign{W_j~&=~\prod_{t=1}^g(f_{2t-1}f_{2t} - b_j)\cr
& =~ (-1)^g\left[b^g_j -
\left( \sum_1^gf_{2t-1}f_{2t}\right) b^{g-1}_j
 + \cdots + (-1)^g\left(\prod_1^{2g}f_t\right)\right]\cr}
\leqno{7.7}$$
which generates the kernel of the restriction map to $H^*(E_{i-1,i};{\bf A})$
and $H^*(E_{i,i-1};{\bf A})$ respectively according to Macdonald's theorem
(6.4).

\noindent{\bf Lemma} 7.8:
\tensl In the stable range the differential $d_1$ in the
spectral sequence (4.3) vanishes for $n>1$.
For $n=1$, it has the form
$$d_1(W_0W_1b_0^{i-2g}b_1^{i-2g}) =
2\left(\sum_1^g f_{2i-1}f_{2i}\right)(W_0W_1b_0^{i-2g-1}b_1^{i-2g-1})
|M_g^*|$$
and it is multiplicative in the remaining generators. That is to say,
$$\Theta W_0W_1b_0^{i-2g}b_1^{i-2g} \fract{d_1}{\ra{1.5}}
\pm \Theta d_1(W_0W_1b_0^{i-2g}b_1^{i-2g}).$$
\tenrm

\noindent{\sc Proof}:  We check the action map $M\times E_{i-1, \dots, i-1} 
\ra{1}E_{i,\dots, i}$ in cohomology.  The term $f_{2r-1}f_{2r} - b_i$ maps to
$$ f_{2r-1}\otimes f_{2r} - f_{2r}\otimes f_{2r-1} + 
1\otimes (f_{2r-1}f_{2r} - b_i)$$
since $f_{2r-1}f_{2r} - b = 0$ in $H^2(M_g;{\bf Z})$.  Hence
$$W_i \mapsto \sum_{r=1}^g (f_{2r-1}\otimes f_{2r} - f_{2r}\otimes f_{2r-1})1
\otimes W_i(r)+ 1\otimes W_i$$
where, as is evident $W_i(r) = \prod _{j \ne r} (f_{2j-1}f_{2j} - b_i)$.  
Thus, since $f_i f_j = 0$ in $H^*(M_g;{\bf Z})$ for $\langle i, j\rangle \ne 
\langle 2r-1, 2r\rangle$ for some $r$, it follows that
$$\eqalign{W_0W_1 ~&\mapsto~ \sum_{r=1}^g 2b \otimes
\left(f_{2r-1}f_{2r} W_0(r)W_1(r)\right) \cr
&\hbox to .3in{\hfill} +
(f_{2r-1}\otimes f_{2r} - f_{2r}\otimes f_{2r-1})1\otimes (W_0(r)W_1 + 
W_0W_1(r))\cr
&\hbox to .3in{\hfill} + 1\otimes W_0W_1.\cr}$$
Next note that $W_i b_i^{i-2g} = 0$ in $H^*(E_{i-1,\dots,i-1}; {\bf
Z})$, so the image
of $(b_0b_1)^{i-2g}W_0W_1$ is $2b\otimes \sum f_{2r-1}f_{2r} b_0^{i-2g-1}
b_1^{i-2g-1}
W_0W_1$ since $b_iW_i(r)f_{2r-1}f_{2r} = W_if_{2r-1}f_{2r}$.  This proves the
desired result.\hfill\za

\noindent{\bf Remark} 7.9: Actually, what the argument above really determines
is the coaction map in cohomology. Given the action map 
$$SP^r(M_g)\times LE_i\ra{1.5} LE_{i+r}$$
described earlier, and suppose $i\ge 2g$, then
$$U_{i+k} \mapsto b^k\otimes {2^k\over
k!}\left(\sum_1^g f_{2r-1}f_{2r}\right)^kU_i$$
where $U_j$ is the generating class above for $H^{2(j-g)}(LE_j;{\bf Z})$
when $n=1$.


\noindent{\bf\twelverm\S8.  Stabilization}

For $d \ge g$, any $n\geq 1 $ and any neighborhood $U_{\epsilon}(*)$,
there are $n+1$ 
divisors $D_i\subset SP^d(U_{\epsilon}(*))$, $0\le i \le n$ with
$\mu(D_0) = \mu(D_1) = 
\cdots = \mu(D_n)$ and $\bigcap D_i = \emptyset$.  By deforming $M_g-*$ to
$M_g-U_{\epsilon}$ and then adding in the corresponding divisors, one obtains a
stablization map $\tau$ which one can iterate
$$\hol{k}^*(M_g, {\bf P}^n) \fract{\tau}{\ra{1.5}} \hol{k+d}^*(M_g, {\bf P}^n) 
\ra{1.5} \cdots \ra{1,5}\hol{k+rd}^*(M_g, {\bf P}^n)\leqno{8.1}$$
obtaining in the limit a space $\lim_{r\mapsto
\infty}\hol{k+rd}^*(M_g, {\bf P}^n)$  
which is homotopy equivalent to any component in
the mapping space $\map{*}^*(M_g, {\bf P}^n)$. Note that this last statement
is a consequence of Segal's stabilization theorem, [S].

The inclusion $\tau$ descends to a map of spectral sequences and by an
argument similar 
to one in [C$^2$M$^2$], we see the corresponding map on the $E_1$-term
is given by 
$$\eqalign{\cup(b_0b_1\cdots b_n)\otimes& id \colon H^*(LE_i;{\bf Z})\otimes 
H^*(SP^j(\Sigma M_g), SP^{j-1}(\Sigma M_g); {\bf Z})\cr
& \ra{1.5} H^*(LE_{i+1};{\bf Z})\otimes H^*(SP^j(\Sigma M_g),
SP^{j-1}(\Sigma M_g); {\bf Z}).\cr}
\leqno{8.2}$$
This relation between inclusion at the level of mapping spaces and
cupping with $b_0\cdots 
b_n$ in the Poincar\'e duals is quite important and follows basically
by checking normal 
bundles (cf [K1]). In any case the content of the preceeding remarks can be 
summarized in

\noindent{\bf Lemma}
 8.3:  \tensl In the stable range the cohomology of $\hol{k}^*(M_g, {\bf P}^n)$
is isomorphic to the cohomology of $\map{k}^*(M_g, {\bf P}^n)$ via the
natural inclusion. 
\tenrm

\noindent This is again direct from the stabilization via cupping with 
$b_0\cdots b_n$ and Segal's stability theorem.

The differentials which we know appear come from comparison with the 
$Div$-space spectral sequence.  These are the differential $d_1$
already described and the further differentials which only hold $\bmod{(p)}$ 
for odd $p$,
$$\eqalign{d(|\gamma_{p^i}|^*)~& =~ {1\over p^i}
\left(\sum_1^g f_{2j-1}f_{2j}\right)^{p^i},\cr
d(|\gamma_{p^i}^{p-1}|\gamma_{p^i}|^*)~&= ~\left[{1\over p^i}\left(\sum_1^gf_{2j-1}f_{2j}\right)^{(p-1)p^i}\right]|\gamma_{p^i}|^*.\cr}
\leqno{8.4}$$
The last differential in (8.4) is called the Kudo differential, [Ku].
There are no differentials for $p = 2$ (see next remark).
Similarly, for $n>2$ all the differentials are zero for all $p$ and hence, even
over the integers. 

\noindent{\bf Remark} 8.5:  There is a sequence of Serre fibrations
$$\Omega^2({\bf P}^n) \ra{1.5} \map{}^*(M_g, {\bf P}^n) \ra{1.5}
\left(\Omega ({\bf P}^n) 
\right)^{2g}\fract{f^{\dag}}{\ra{1.5}} \Omega ({\bf P}^n)$$
induced from the cofibration sequence associated to $M_g$;
$$S^1\fract{f}{\lrar}\bigvee^{2g}S^1\hookrightarrow M_g\lrar S^2$$
which gives $\map{}^*(M_g, {\bf P}^n)$ as the total space of a principal
$\Omega^2({\bf P}^n)$ 
fibration with classifying map $f^{\dag}$.  Then the discussion above
shows by dualizing 
that in the Serre spectral sequence of the fibration
$$\Omega^2({\bf P}^n) \ra{1.5} \map{}^*(M_g, {\bf P}^n) \ra{1.5}
\left(\Omega ({\bf P}^n)\right)^{2g} 
\leqno{8.6}$$
we have (for $n=1$) 
the differentials $d^1(|b|^*) = 2\sum_1^g f_{2j-1}f_{2j}$, the
transgressive $\bmod{(p)}$ 
differentials from $|b^{p^i}|^*$ to the divided power of
$\left(\sum_1^gf_{2j-1}f_{2j}\right)$, 
and the Kudo differentials for $p$ odd.  Of course, for $p=2$ the
differentials in the Serre 
spectral sequence are totally trangressive from the fiber to the base.
Hence, there are no differentials and
$E_2 = E_{\infty}$.  Similarly, for $n>2$ there are no differentials
for any $p$ and 
$E_2 = E_{\infty}$ in this spectral sequence even with integer
coefficients though the 
only way we know to prove this is via the results of [K2].


\noindent{\bf\twelverm\S9.  The Genus One Case}

When $g=1$, then $\mu\colon M_1 \ra{1} J(M_1)$ is a holomorphic
homeomorphism identifying the Jacobi variety with $T$ itself.
Also, the stable range starts immediately in this case and hence,
for all $n\geq 1$ the Abel-Jacobi map $\mu$ is a fibration 
$${\bf P}^{n-1}\lrar \sp{n}(T)\fract{\mu}{\lrar } T.$$
In the spectral sequence for
$QE_k(n)$ with $n\ge 2$ (dual to $Hol^*_k(T, {\bf P}^n)$) there are
no differentials except for those at the tail end (6.15) and
(6.16) but these 
are the only differentials and $E^1= E^{\infty}$ away from this region
while $E^2 = E^{\infty}$ for the entire spectral sequence.  In the case $n=1$
the only differential is again $d_1$ which now also has a stable component
given in (5.9), and again $E^2 = E^{\infty}$.

It also follows that duality gives the isomorphism
$$ {\tilde H}^{2k(n+1)-2n-*}(\hol{k}^*(T,{\bf P}^n); {\bf F} )\cong
H_*(QE_{\underbrace{k,\ldots,k}_{n+1}};{\bf F} )
\leqno{9.1}$$
for all $k\ge 2$, with $Hol_1^*(T,{\bf P}^n)$ being empty for all $n$.

We now determine these homology groups explicitly for $n=1$.  First,
in the case where ${\bf F}$ has characteristic 2 there are
no stable differentials and $E^1 = E^{\infty}$ in this range.  Second,
in case ${\bf F} = {\bf Q}$ we have
$$\eqalign{E^1(\hol{k}^*)~&=~\coprod_{i=2}^k\Sigma^{4(i-1)}H_*(T;{\bf Q})\otimes
\left({\bf Q}[|e_1|, |e_2|]_{k-i} \oplus |[T]|\otimes {\bf Q}[|e_1|,
|e_2|]_{k-i-1} \right)\cr
&\hbox to .3in{\hfill}\oplus
{\tilde H}_*(T;{\bf Q})\otimes
\left({\bf Q}[|e_1|, |e_2|]_{k-1} \oplus |[T]|\otimes {\bf Q}[|e_1|,
|e_2|]_{k-2} 
\right)\cr
&\hbox to .6in{\hfill}\oplus {\bf Q}[|e_1|, |e_2|]_k + |T|{\bf
Q}[|e_1|, |e_2|]_{k-1} 
.\cr}
\leqno{9.2}$$
where ${\bf Q}[|e_1|, |e_2|]_r$ is the $(r+1)$-dimensional subspace
spanned by the 
monomials of degree $r$; 
that is: $|e_1|^r, ~|e_1|^{r-1}|e_2|, \dots, |e_2|^r$.  The
only differential is $d_1$, generated by
$$d_1(\Sigma^{4(k-j-1)} e_1e_2 |T|) = \Sigma^{4(k-j)} 1$$
and the unstable differentials of (6.15) and (6.16).
Thus one has directly (compare (6.19) for the last term) that
$$\eqalign{E^2~&=~
\coprod_{i=2}^k\Sigma^{4(i-1)}{\tilde H}_*(T;{\bf Q})\otimes
{\bf Q}[|e_1|, |e_2|]_{k-i}\oplus
\Sigma^{4(i-1)}|T|(1, e_1, e_2){\bf Q}[|e_1|, |e_2|]_{k-i-1}\cr
&\hbox to .3in{\hfill} \oplus {\bf Q}\{ e_2 |e_1|^{k-1}, e_2|e_1|^{k-2}|e_2|, 
\dots,e_2|e_1||e_2|^{k-2}\}\cr
&\hbox to .6in{\hfill}\oplus |T|(1, e_1,e_2){\bf Q}[|e_1|, |e_2|]_{k-2}\cr
&=~E^{\infty}\cr}
\leqno{9.3}$$
In the remaining case where ${\bf F} = {\bf F}_{p^n}$ has characteristic $p$ 
an odd prime,
as explained in (6.6), (6.7) we can write the $E_1$-term as
$$\coprod_{i + j = k}\Sigma^{4(i-1)}H^*(T; {\bf F})\otimes
\sum_{l=0}^j\Sigma^{4l}{\cal D}_j^*(p)\otimes
{\bf F}[|b_1, b_{2}]_{j-l}$$
where ${\cal D}_j^* = 0$ if $j\ne pk$ or $pk+1$ and with the usual
special considerations 
at the tail end of the filtration.  Incidently, the case $j=pk+1$ implies
that $\Sigma^{4j}{\cal D}_j^* = |M|\Sigma^{4(j-1)}{\cal D}_{j-1}^*$.  The
differential 
is as before and we have that in the stable range
$$\eqalign{E^2~& =~ \coprod_{i=0}^k\Sigma^{4(i+j-1)}{\tilde
H}_*(T;{\bf F})\otimes {\cal D}_j^* 
\otimes {\bf F}[b_1,b_2]_{k-i-j} \oplus
\Sigma^{4(i+j-2)}|T|{\bf F}[b_1,b_2]_{k-i-j-1}\cr 
&=E^{\infty}\cr}$$

In each case the homology injects into $H_*(\map{k}^*(T, {\bf
P}^1);{\bf F})$. Thus the geometric interpretation of each of the
homology and cohomology classes in $H_*(\hol{k}^*(T,{\bf P}^1);{\bf
F})$ can be regarded as the same as that for the corresponding class
in the mapping space. 

For example with rational coefficients we have the following corollary
to 9.3 and 9.1

\noindent{\bf Lemma} 9.4:\tensl~
Let $I = (|e_1|, |e_2|)$ be the augmentation ideal in ${\bf Q}[|e_1|, |e_2|]$. 
Then
$$\eqalign{
H_*(\hol{k}^*(T, {\bf P}^1);{\bf Q}) ~&\cong ~{\bf Q}(1,
e_1, e_2, h_{2,1}, h_{2,2}, v_3)\otimes
{\bf Q}[|e_1|, |e_2|]/I^{k-1}\cr
&\hbox to .3in{\hfill} \oplus {\bf Q}\{e_2|e_1|^{k-1}, \dots,
e_2|e_1||e_2|^{k-2}\}\cr} $$
\tenrm

In particular the Poincar\'e series
for $H_*(\hol{k}^*(T, {\bf P}^1);{\bf Q})$ is
$$(1+2x+2x^2+x^3)(1 + 2x^2 + 3x^4 + \cdots + (k-1)x^{2(k-2)}) +
(k-2)x^{2(k-1)-1}.$$ 
Note how this compares with the corresponding Poincar\'e series for
the mapping space 
$H_*(\map{}^*(T, {\bf P}^1);{\bf Q})$:
$${(1+2x+2x^2 +x^3)\over (1-x^2)^2}.$$
{\tensl The case $n\ge 2$.}  Any component in
the mapping space
is given as the total space of the Serre fibration
$$\Omega^2S^{2n+1} \ra{1} E \ra{1} (S^1\times \Omega S^{2n+1})^2$$
with
$$H_*(\map{}^*(T, {\bf P}^n);{\bf F}) = H_*(\Omega^2S^{2n+1};{\bf F})\otimes 
H_*(T;{\bf F})\otimes {\bf F}[g_{2n}, g_{2n}^{\prime}].$$
Consequently, with ${\bf Q}$ as coefficients the Poincar\'e series for
the mapping space has the form
$${(1+x)^2(1+x^{2n-1})\over (1-x^{2n})^2}$$
while for $\hol{k}^*(T, {\bf P}^n)$ we have
$$(1+x)^2(1+x^{2n-1})(1 + 2x^{2n} + 3x^{4n} + \cdots + (k-1)x^{(k-2)2n})
+(k-2)x^{2(k-1)n-1}.$$

\noindent{\bf Remark} 9.5:  Notice that in this case all the
calculations were  
independent of the particular elliptic curve $T$ we were studying, even though 
different $T$ are
not holomorphically isomorphic.  In fact this is not unexpected since we have

\noindent{\bf Theorem} 9.6:
{\sl~Let T be any non-singular elliptic curve, then the homeomorphism type of
$\hol{k}^*(T,{\bf P}^n)$ is independent of $T$.}

\noindent{\sc Proof:}  Let $T$ and $T'$ be two given tori
with corresponding lattices $\Lambda $ and $\Lambda '\subset{\bf C}$.
Then there is a map $\Lambda\lrar\Lambda'$ 
which induces a homeomorphism
$$T={\bf C}/\Lambda\fract{f}{\lrar} T'={\bf C}/\Lambda'.$$
and such that the following commutes for all $n$;
$$\matrix{\sp{n}(T)&\fract{\sp{n}(f)}{\ra 2}&\sp{n}(T')\cr
\decdnar{\mu}&&\decdnar{\mu}\cr
T&\fract{f}{\ra 2}&T'.\cr}$$
It follows therefore that $E_{i,i,\ldots i} \cong E_{i,i,\ldots, i}'$ where
$E_{i,i,\ldots i}$ (resp. $E_{i,i,\ldots, i}'$) is the pushout as in \S 3
$$\matrix{E_{i,i,\ldots,i}&\lrar&\sp{i}(T)\times\cdots\times\sp{i}(T)\cr
\downarrow&&\downarrow\cr
T&\fract{\Delta}{\lrar}&T\times\cdots\times T\cr}$$
(respectively, $T$ is replaced by $T'$). It is now clear that
$$\hol{n}^*(T,{\bf P}^n) \cong {\tilde E}_{i,\ldots i}-\hbox{Image}(\nu )
\cong {\tilde E}'_{i,\ldots i}-\hbox{Image}(\nu) \cong 
\hol{n}^*(T,{\bf P}^n)$$
and the proof is complete.
\hfill\za


\noindent{\bf\twelverm\S10. The Decomposition for Hyperelliptic Curves}

The surface $M_g$ is hyperelliptic if and only if there is a degree two 
(branched) holomorphic map
$f\colon M_g \ra{1} {\bf P}^1$. This map $f$ fits in the following diagram
$$\matrix{M_g&&\fract{s}{\ra 3}&&{\bf P}^{2g-2}\cr
&\searrow&&\nearrow\cr
&&{\bf P}^1&\cr}\leqno{10.1}$$
where $s$ is the canonical map, [GH, p. 247]. 

From the point of view of the Abel-Jacobi maps, $M_g$ is hyperelliptic
if and only if 
we have a cofibering of the form
$${\bf P}^1 \hookrightarrow \sp{2}(M_g) \ra{1.5} \sp{2}(M_g)/{\bf P}^1
= W_2\subset J(M_g).$$ 
The holomorphic embedding ${\bf P}^1 \hookrightarrow SP^2(M_g)$ is now
constructed 
by associating to $z\in {\bf P}^1$ the pair $f^{-1}(z) \in SP^2(M_g)$
where the points are counted with multiplicity, i.e. the ramification
points are counted twice.  The image $\tau \in J(M_g)$ of ${\bf P}^1$ under
$\mu_2$ is called the hyperelliptic point of $M_g$.

\noindent{\bf Lemma}  10.2:  \tensl  Suppose that $M_g$ is 
hyperelliptic and $\tau \in J(M_g)$ is the hyperelliptic point.
Let $k\geq 1$ and $t\leq \left[ {k\over 2}\right]$. Then
\item{(a)}  for $k \le g$, the space $W_k^t$
is $\mu_{k-2t}(SP^{k-2t}(M_g)) + t\tau$,
\item{(b)}  for $2g-1 > k > g$, $t>k-g$, we also have
$W_k^t=\mu_{k-2t}(SP^{k-2t}(M_g)) + t\tau$.
\tenrm

\noindent{\sc Proof}:  Note that $\sp{r}({\bf P}^1) = {\bf P}^r$ so the pairing
$$\sp{2l}(M_g)\times \sp{k-2l}(M_g) \ra{1.5} \sp{k}(M_g)$$
induces an action
$$\nu\colon~
{\bf P}^l\times \sp{k-2l}(M_g) \fract{=}{\ra{1.5}} \sp{l}({\bf
P}^1)\times \sp{k-2l}(M_g) \ra{1.5} \sp{k}(M_g)$$
and the image of
$$\nu(\langle p_1, \dots, p_l\rangle, \langle m_1, \dots, m_{k-2l}\rangle)$$
under the Abel-Jacobi map is $l\tau + \mu_{k-2l}(\langle m_1, \dots,
m_{k-2l}\rangle$. 
From this the fact that $W_k^t$ is at least a large as asserted in (a)
follows directly. 
The inclusion argument for (b) is a similar dimension check when we recall
that generically $\mu_{g+s}^{-1} (j) = {\bf P}^s$, so we are checking
for points with inverse image ${\bf P}^t$ where $t> s$
and $\nu({\bf P}^{t}\times \sp{g+s -2t}(M_g)) 
\subset \mu^{-1}_{g+s}(W_{g+s}^{t})$.

The converse is given on page 13 of [ACGH] for $k \le g$ where it is
pointed out that as 
a consequence of the geometric version of the Riemann-Roch theorem
every ${\bf P}^t$ 
inverse image under the Abel-Jacobi map for $\sp{k}(M_g)$ has
the form $\mu^{-1}(t\tau) + m_1 + \cdots + m_{k-2t}$ where no two of
the $m_i$'s 
are conjugate in $M_g$ under the hyperelliptic involution.  To obtain
the remaining 
part of the converse one uses the residual series isomorphism
$$W_k^t\cong W^{g-k+t-1}_{2g-2-k}, \hbox to .3in{\hfill} D
\leftrightarrow K - D$$ 
which holds for all curves $M_g$, $k\geq g$.\hfill\za

For later use we also specify an appropriate base point for the construction.  
Let $*\in M_g$ be one of
the ramification points of $f\colon M_g \ra{1} {\bf P}^1$,
then $2\mu(*)= \tau$,
and we choose $*$ as the base point in $M_g$.  Also the base point of
${\bf P}^1$ is $\langle *,*\rangle$ and the inclusion $f\colon {\bf P}^1 \ra{1}
\sp{2}(M_g)$ is based.  With respect to this choice, the
composite inclusion
$$\sp{k}(M_g) \hookrightarrow \sp{k+1}(M_g)\hookrightarrow \sp{k+2}(M_g)$$
commutes with the ${\bf P}^1$-action.  Moreover, in the actual
construction of the 
Abel-Jacobi map $\mu$ given at the beginning of \S1~
we can take our base point in $M_g$ to be $*$, so
that $\mu(*) = 0 \in J(M_g)$, and the hyperelliptic point $\tau = 0$
as well.  This 
has the advantage that the inclusion $W_k^t + \tau \subset
W_{k+2}^{t+1}$ becomes 
simply $W_k^t \subset W_{k+1}^{t+1}$.


\noindent{\bf\twelverm \S11.  Some preliminary constructions}

Let $V$ be a space with a commutative, associative multiplication.  For example
$$V = SP^{\infty}(X,*),~ V = \coprod_1^{\infty}\sp{n}(X),~ V =
\coprod_1^{\infty}  
\sp{i}(X)\times \sp{i}(X),$$
etc.  If $f\colon X \ra{1} V$ is any continuous map
then $f$ defines an action
$$\mu_f\colon \coprod \sp{k}(X)\times V \ra{1.5} V$$
by $\mu_f(\langle x_1,\dots,x_k\rangle, v) = f(x_1)f(x_2)\cdots f(x_k)v$.

Suppose now that $V$, $f_1\colon {\bf P}^1 \ra{1} V$ and $V^{\prime}$, $f_2\colon
{\bf P}^1 \ra{1} V^{\prime}$ are given with $V$, $V^{\prime}$ as
above.  Then we can form the bicomplex
$$B(f_1,f_2)~\colon = ~V\times_{f_1}\left(\coprod_k \sp{k}(I\times
{\bf P}^1)\right) 
\times_{f_2}V^{\prime},\leqno{11.1}$$
defined by the equivalence relation built from the basic relations
$$\eqalign{(v,&\langle(t_1,p_1),\dots, (t_k,p_k)\rangle ,v^{\prime})\sim \cr
&\hbox to .3in{\hfill}\cases{(v_1\cdot f_1(p_i),
\langle (t_1,p_1),\dots, \widehat{(t_i,p_i)}, \dots, (t_k,p_k)\rangle,
v^{\prime})&  if $t_i = 0$,\cr
(v_1,\langle (t_1,p_1),\dots, \widehat{(t_i,p_i)}, \dots,
(t_k,p_k)\rangle, f_2(p_i)\cdot v^{\prime})& if $t_i = 1$.\cr}\cr}$$

\noindent{\bf Remark} 11.2:  This is just the 2-sided bar construction
which is associated 
to the homological functor $Tor_{A}(M_A, _A N)$.  But here,
left and right modules are equivalent since $A$ is commutative.

\noindent{\bf Remark} 11.3:  This is an unbased construction. If $f_1$
and $f_2$ are based 
maps with $f_i(*) = id$, $i = 1,2$, then there is an associated
reduced construction, 
${\tilde B}(f_1, f_2)$.  Also, the fact that we are using ${\bf P}^1$
plays no real role 
in the definition.  It is just there because that is all we
use in the applications of this construction to our study of
$Hol_k^*(M_g, {\bf P}^n)$ for $M_g$ hyperelliptic.

\noindent{\bf Example} 11.4:  Suppose that $V = \coprod \sp{k}(X)$ and
$f_1\colon 
{\bf P}^1 \ra{1} \sp{i}(X)$, then write
$$V_j ~=~\coprod_{k=0}^{\infty}\sp{j+ik}(X),\hbox to .3in{\hfill} 0\le j<i$$
so we have separate actions on each $V_j$.

\noindent{\bf Example} 11.5:  Suppose $f_2\colon {\bf P}^1 \ra{1} *$
where $*$ has the obvious commutative and associative action.  Then
$$B(f_1,f_2) = V\times_{f_1}\coprod\sp{n}(I\times {\bf P}^1)\times_{f_2}*
\simeq V\times_{f_1}\coprod\sp{n}(c{\bf P}^1)$$
is a space we have seen many times before.

\noindent{\bf Example} 11.6:  $\Delta_{n+1}\colon {\bf P}^1 \ra{1.5}
({\bf P}^1)^{n+1}$ induces 
an action on $\left(\coprod \sp{k}({\bf P}^1)\right)^{n+1}$ which
restricts, in a manner 
similar to that in example (1) to an action on
$${\cal E}_n(i_1, \dots, i_{n+1}) ~=~
\coprod_{k=0}^{\infty}\sp{i_1+k}({\bf P}^1)\times 
\cdots \times \sp{i_{n+1} + k}({\bf P}^1)$$
for each $(n+1)$-tuple $(i_1,\dots , i_{n+1})$ of non-negative
integers with at least one of the $i_j = 0$.


\noindent{\bf\twelverm \S12.  Models for the spaces
$E_{i,\dots, i}$, $LE_{i,\dots, i}$ and $W_k^t$ when $M_g$ is
hyperelliptic}

We specialize to
$$BH(f,n)~=~\left(\coprod_m\sp{m}(M_g)\right)\times_{f}
\left(\coprod_m\sp{m}(I\times
{\bf P}^1)\right)\times_{f_2}{\cal E}_n(0,0,\dots, 0)\leqno{12.1}$$
where $f\colon {\bf P}^1 \hookrightarrow \sp{2}(M_g)$ embeds ${\bf
P}^1$ as the inverse 
image under the Abel-Jacobi map of the hyperelliptic point, $0$, for
some hyperelliptic 
structure on $M_g$ and $f_2$ is the structure map $\Delta_{n+1}$ of
(11.6).

This construction specializes, as in (11.5), to give two disjoint subspaces,
$$\eqalign{12.2\hbox to .5in{\hfill}
BH_0(f,n)~&=~V_0\times_f \left(\coprod_m\sp{m}(I\times {\bf P}^1)\right)
\times_{f_2}{\cal E}_n(0,\dots, 0)\cr
12.3\hbox to .5in{\hfill}
BH_1(f,n)~&=~V_1\times_f\left(\coprod_m\sp{m}(I\times {\bf
P}^1)\right)\times_{f_2} {\cal E}_n(0,\dots, 0)\cr}$$
where $V_0 = \coprod \sp{2k}(M_g)$ and $V_1 = \coprod\sp{2k+1}(M_g)$.  
Each $BH_{\epsilon}(f,n)$ breaks up into further components as follows.
We define a grading by letting
$$\eqalign{(\langle m_1, \dots, m_l\rangle, \langle (t_1,p_1), &\dots,
(t_r,p_r)\rangle, 
(w_{s,1},\dots, w_{s,n+1}))\cr
&\in \sp{l}(M_g)\times \sp{r}(I\times {\bf P}^1)\times
\prod_1^{n+1}\sp{s}({\bf P}^1)\cr}\leqno{12.4}$$
have grading degree $l + 2r + 2s$.  This grading is obviously preserved by the
identifications and hence induces the desired decomposition of the
quotient $BH_{\epsilon}(f)$. 
We write ${\cal G}_s(f,n)$ for the component of points of grading degree $s$.

This is the general case, however, there is a special case
$BH_{\epsilon}(f,-1)$ which is 
also needed where $V_2 \sim *$ consists of a single point.  In this
case the grading 
degree is {\tensl not preserved} by the identifications, and the
grading only gives a 
{\tensl filtration} of $BH_{\epsilon}(f,-1)$.  We write ${\cal
F}_k(BH_{\epsilon}(f,-1))$ 
for the filtering subspace consisting of all points of grading degree
$\le k$.  Moreover, 
in this case the identification has the form $(v, (t, x), *)\sim (vf(x), *)$
when $t = 0$ and $(v, (t,x), *) \sim (v,*)$ when $t = 1$.

\noindent{\bf Lemma} 12.5:  \tensl We have the following homotopy
equivalences in the 
case where $M_g$ is hyperelliptic.
\item{(a)} ${\cal G}_k(BH_{\epsilon}(f,n)) \simeq E_{k,\dots, k}$ for
$k\equiv \epsilon\bmod{(2)}$ provided $k \le g$.
\item{(b)} ${\cal F}_k(BH_{\epsilon}(f,-1)) \simeq W_k$, $k\equiv
\epsilon\bmod{(2)}$ and $k< g$.
\tenrm

\noindent(These constructions are modeled on the analogous constructions in
[M2] and the proofs are the same.)

\noindent{\bf Remark} 12.6:  The point of these constructions is to
handle the algebraic 
complexities of systematically collapsing out ${\bf P}^n$'s.  As one
sees from the 
lemma, the effect is to introduce $\sp{n}(\Sigma {\bf P}^1)$'s which
provide a measure 
of the geometry of the collapsed spaces as the collapsing gets more and more
complex.


\noindent{\bf\twelverm \S13.  Some spectral sequences for the
$BH_{\epsilon}(f,n)$} 

We can bifilter ${\cal G}_k(BH_{\epsilon}(f,n))$ and ${\cal
F}_k(BH_{\epsilon}(f,-1))$
by the number of $\sp{t}(I\times {\bf P}^1)$
terms appearing, 
obtaining spectral sequences converging to the various spaces above as
well as certain useful quotients.  The $E^1$-terms are as follows:
$$\coprod_{k+2l+2t= v}H_*(\sp{k}(M_g))
\otimes H_*(\sp{l}(\Sigma {\bf P}^1_+), \sp{l-1}(\Sigma {\bf P}^1_+))\otimes
\left(H_*({\bf P}^t)\right)^{n+1}\leqno{13.1}$$
converging to $H_*(E_{v,\dots, v})$ with field coefficients and $v<g$,
$$\coprod_{k + 2l \le v}H_*(\sp{k}(M_g))\otimes
H_*(\sp{l}(\Sigma {\bf P}^1_+), \sp{l-1}(\Sigma {\bf P}^1_+))\leqno{13.2}$$
with $k \equiv v\bmod{(2)}$ converging to $H_*(W_v)$ for $v<g$, and
$$\coprod_{k+2l = v} H_*(\sp{k}(M_g))\otimes
H_*(\sp{l}(\Sigma {\bf P}^1_+), \sp{l-1}(\Sigma {\bf P}^1_+))\leqno{13.3}$$
converging to $H_*(W_v/W_v^1)$.  (The notation $\Sigma {\bf P}^1_+$
means the reduced suspension on the union of ${\bf P}^1$ with a disjoint
base point denoted $+$, it has the homotopy type of the wedge
$\Sigma({\bf P}^1) \vee S^1$.)

If we relativize we obtain a spectral sequence with $E^1$-term
$$\eqalign{\coprod_{k + 2(l+t) = v}H_*(\sp{k}(M_g),&\sp{k-1}(M_g))
\otimes H_*(\sp{l}(\Sigma {\bf P}^1_+),\sp{l-1}({\bf P}^1_+))\otimes
{\tilde H}_*(S^{2(n+1)t})\cr
&\Longrightarrow {\tilde H}_*(LE_{v,\dots, v}).\cr}\leqno{13.4}$$
All these spectral sequences are algebraic in the sense of [K1], and can be
modeled by simply replacing the terms $H_*(\sp{l}(\Sigma {\bf P}^1))$ by
the cobar construction on
$$H_*\left(\coprod \sp{l}({\bf P}^1_+)\right) = {\bf Z}[t]\otimes \power(b).$$
In particular, the relative spectral sequence splits as a direct sum
of spectral sequences each of the form
$$\eqalign{&\left[\coprod_{k+2l = v - 2t}H_*(\sp{k}(M_g),\sp{k-1}(M_g))
\otimes H_*(\sp{l}(\Sigma {\bf P}^1_+),\sp{l-1}(\Sigma {\bf P}^1_+))\right]
\otimes {\tilde H}_*(S^{2t(n+1)})\cr
& \hbox to 2in{\hfill}\Longrightarrow
\Sigma^{2t(n+1)}\left(H_*(W_{v-2t}/(W^1_{v-2t}\cup W_{v-2t-1}))\right)\cr}
\leqno{13.5}$$
for any $n >0$.

In the case where the coefficients are the rationals
these spectral sequences becomes quite simple.  First
$H_*(\sp{\infty}(\Sigma {\bf P}^1_+);{\bf Q})= \ext(h_1,h_3)$, and 
$$H_*(\sp{i}(\Sigma {\bf P}^1), \sp{i-1}(\Sigma {\bf P}^1);{\bf Q}) = 0$$
for $i \ge 3$.  Second, when we write
$$H_*(\sp{k}(M_g), \sp{k-1}(M_g);{\bf Z}) =
\bigoplus_{s=0}^k\ext_s(e_1, \dots, e_{2g})\otimes 
[M_g]^{k-s},$$
for $k < g$, we see that the action map
$$*[{\bf P}^1]\colon H_*(\sp{k}(M_g), \sp{k-1}(M_g);
{\bf Z}) \ra{1} H_{*+2}(\sp{k+2}(M_g), \sp{k-1}(M_g);{\bf Z})$$
is just multiplication by $-\sum_{j=1}^ge_{2j-1}e_{2j}$ since
$$[{\bf P}^1] = [M_g] -
\sum_1^ge_{2j-1}e_{2j}$$
in $H_*(\sp{2}(M_g);{\bf Z})$.  Third, from [BCM], this map is
injective in homology for $k \le g-2$.  Thus we have

\noindent{\bf Lemma} 13.6: \tensl 
$$\eqalign{H_*(LE_{v,\dots, v}&;{\bf Q}) =\cr
& \sum_{t=0}^{[v/2]} \Sigma^{2t(n+1)}H_*(\sp{v-2t}(M_g),
\sp{v-2t-1}(M_g);{\bf Q})/im\left(*(\sum_1^ge_{2j-1}e_{2j})\right)\cr}$$
for $v\le g$, which has Poincar\'e series
$$\sum_{t=0}^{[v/2]}x^{2t(n+1)}S_{v-2t}$$
where $S_r = \sum_{l \le r}\left[{2g\choose l} - {2g\choose
l-1}\right]x^{2r-l}$.\tenrm 

Similarly, we can analyze $H_*(W_j;{\bf Q})$.  We have

\noindent{\bf Lemma} 13.7:  \tensl
\item{(a)}  The inclusion $W_j \subset J(M_g)$ induces an injection in rational
homology $H_*(W_j;{\bf Q}) \hookrightarrow H_*(J(M_g);{\bf Q}) =
\ext(e_1, \dots, e_{2g})$ 
with image the subvector space spanned by the subspaces
$$\left\{ \ext_s(e_1,\dots,e_{2g})[M_g]^t~|~s + t \le j\right\}$$
where $[M_g] = \sum_1^g e_{2i-1}e_{2i}$ is the image of the
fundamental class of $M_g$ 
under the Abel-Jacobi map $\mu_*$,
\item{(b)}  $H_*(W_{j-1};{\bf Q})$ injects into $H_*(W_j;{\bf Q})$
under the inclusion so the relative groups are given as
$$H_*(W_j, W_{j-1};{\bf Q}) \cong H_*(W_j;{\bf Q})/H_*(W_{j-1};{\bf Q}).$$
\tenrm

\noindent{\sc Proof}: 
We have that the $E^1$-term of the spectral sequence above is given as
$$\eqalign{H_*(\sp{j}(M_g);{\bf Q}) &\oplus H_*(\sp{j-2}(M_g);{\bf
Q})(1, h_1, h_3)\cr 
&\oplus \sum_{2\le 2l \le j}H_*(\sp{j-2l}(M_g);{\bf Q})(1, h_1, h_3,
h_1h_3).\cr}$$ 
$d_1(\Theta h_1) = i_*(\Theta) - \Theta$ where $\Theta \in H_*(\sp{j-2l}(M_g))$
and $i\colon \sp{j-2l}(M_g) \ra{1} \sp{j-2l+2}(M_g)$ is a base point inclusion.
Similarly $d_1(\Theta h_1 h_3) = d_1(\Theta h_1)h_3 \pm (\Theta h_1)d_1(h_3)$
where, of course, this last term is just multiplication by $[M_g] -
\sum e_{2i-1}e_{2i}$. 
Ignoring the last piece one sees that the $d_1$-differential on $h_1$
reduces the 
calculation of $E^2$ to the following complex:
$$\eqalign{H_*(\sp{j-2}(M_g);{\bf Q})h_3 \ra{1.5} H_*&(\sp{j}(M_g);{\bf Q}),\cr
&\Theta h_3 \mapsto \Theta ([M_g] - \sum_1^g e_{2i-1}e_{2i})\cr}$$
from which the first statement of the lemma follows.  The second
statement is then 
immediate.\hfill\za

\noindent{\tensl The remaining groups $H_*(LE_{v,\dots,v};{\bf Q})$ in the
hyperelliptic case}
\medskip
It is convenient to write $v = g + s$, $1 \le s < g-1$ when $v>g$.  In this
case we have that $J(M_g)$ is filtered by $W_0 \subset W_1 \subset
\cdots W_{g-1} \subset 
J(M_g)$, and over points in $W_{g-s -2l}  - W_{g-s - 2l-2}$ the inverse image
in $E_{v,\dots, v}$ is $({\bf P}^{s+l})^{n+1}$, while the inverse image
over $J(M_g) - W_{g-s-2}$ 
is $({\bf P}^s)^{n+1}$.  Thus, we obtain a spectral sequence with $E^1$-term
$$\eqalign{H_*(J(M_g), W_{g-s-2})\otimes &\left(H_*({\bf P}^s)\right)^{n+1}\cr
& \oplus \sum_{l=1}^{[g-s/2]} H_*(W_{g-s-2l},
W_{g-s-2l-2})\otimes \left(H_*({\bf P}^{s+l})\right)^{n+1}\cr
&\hbox to 1in{\hfill} \Longrightarrow H_*(E_{g+s, \dots,
g+s})\cr}\leqno{13.8}$$
with field coefficients.

Similarly, the corresponding $E^1$-term for $H_*(LE_{g+s, \dots,
g+s})$ is given by 
$$\eqalign{H_*(J(M_g), W_{g-s-1})\otimes &\left(H_*({\bf P}^s, {\bf
P}^{s-1})\right)^{n+1}\cr 
&\oplus \sum_{l=1}^{[g-s/2]} H_*(W_{g-s-2l}, W_{g-s-2l-1})\otimes
\left(H_*({\bf P}^{s+l}, 
{\bf P}^{s+l-1})\right)^{n+1}\cr
&\hbox to 1in{\hfill} \Longrightarrow H_*(LE_{g+s, \dots, g+s}).\cr}
\leqno{13.9}$$

Note here that in the formula above the term involving the relative groups
$$H_*(W_{g-s-2l-1},W_{g-s-2l-2})$$
is zero.  That is because, in this region of $J(M_g)$ and this range
for $v$ the inverse image of a point is ${\bf P}^{s+l}$ in both
$E_{g+s, \dots, g+s}$ and 
$E_{g+s-1, \dots, g+s-1}$.

Again, because of this gap the spectral sequence collapses and
$E^1 = E^{\infty}$.  Of course the calculations become considerably
more complex with ${\bf F}_p$ coefficients.


\noindent{\bf\twelverm\S14.  The $Ext$-groups for $\ext_g/S_g$}

By inspection, in the spectral sequence of (5.9) we find the that
$d_1$-differential gives rise to a series of complexes that occur
in the calculation of $Ext$-groups for the ring $V_g$ defined below.  
In this section we calculate these $Ext$-groups and apply the
results in \S15~to determine the structure of the rational homology
of the space $Hol_k^*(M_g, {\bf P}^n)$ with $M_g$ hyperelliptic and
$n\ge 2$.

\noindent{\bf Notation} 14.1:
\item{(1)} $\ext_g = \ext(e_1, \dots, e_{2g}) = H_*(J(M_g), {\bf Q})$,
\item{(2)} $f_g = \sum_{i=1}^ge_{2i-1}e_{2i}\in \ext_g$,
\item{(3)} $S_g\subset \ext_g$ is the ideal $(f_g)$, and
\item{(4)} $\ext_g/S_g = V_g$ is the quotient.

$V_g$ occurs quite often in both our study of $\hol{d}^*(M_g, {\bf P}^n)$ and
$\map{d}^*(M_g,{\bf P}^n)$.  In particular we need the groups
$Ext_{\ext_g}^{*,*}(V_g, {\bf Q})$ which we will describe as modules over
$Ext_{\ext_g}^{*,*}({\bf Q},{\bf Q}) = {\bf Q}[h_1, \dots, h_{2g}]$ where $h_i
= |e_i|^*$, $1\le i\le 2g$.

In [BCM] is is shown that
$$* f_g\colon (\ext_g)_i \ra{1.5} (\ext_g)_{i+2}$$
is injective for $i \le g-1$ and surjective for $i\ge g-1$.  Consequently,
$(V_g)_i = 0$ for $i>g$ while $(V_g)_i = {\bf Q}^{N(g,i)}$ for $0\le i<g$ where
$$N(g,i)~=~\cases{{2g\choose i} - {2g\choose i-2}& for $2 \le i \le g$,\cr
2g& for $i = 1$,\cr
1& for $i=0$.\cr}$$
As a special case note that
$$V_1 ~=~ \ext_1/\{e_1e_2\},$$
so there is a short exact sequence of $\ext_1$-modules
$$0 \ra{1.5} {\bf Q} \fract{e_1e_2}{\hookrightarrow} \ext_1 \ra{1.5} V_1
\ra{1.5} 0$$ 
where ${\bf Q}$ denotes the ``trivial" module given by the augmentation
map $\epsilon\colon 
\ext_g \ra{1} {\bf Q}$, $\epsilon(e_i) = 0$, $\epsilon(1) = 1$.  This exact
sequence gives rise 
to the long exact sequence of $Ext$-groups:
$$\cdots \ra{1.5} Ext_{\ext_1}^{i-1, j-2}({\bf Q},{\bf Q})
\fract{\delta}{\ra{1.5}} 
Ext_{\ext_1}^{i,j}(V_1, {\bf Q}) \ra{1.5} Ext_{\ext_1}^{i,j}(\ext_1, {\bf Q})
\ra{1.5} Ext_{\ext_1}^{i,j}({\bf Q}, {\bf Q}) \fract{\delta}{\ra{1.5}} \cdots$$
and since $Ext_{\ext_1}^{i,j}(\ext_1,{\bf Q}) = \cases{0& $(i,j) \ne (0,0)$,\cr
{\bf Q}& $(i,j) = (0,0)$,\cr}$ we have
$$\eqalign{Ext_{\ext_1}^{*,*}(V_1,{\bf Q}) ~&=~ {\bf Q} e_{0,0} \oplus
Ext_{\ext_1}^{*,*}({\bf Q}, {\bf Q})P_{1,2}\cr
&=~{\bf Q} e_{0,0} \oplus {\bf Q}[h_1, h_2]P_{1,2}\cr}$$
where $P_{1,2}$ is a generating element in $Ext_{\ext_1}^{1,2}(V_1, {\bf Q})$.

Before we can describe the structure of the modules
$Ext_{\ext_g}^{*,*}(V_g, {\bf Q})$ for $g>1$ we need some preliminary
constructions. 

First, we will often be dealing with a situation where we have a ring $R_{i+1}
= R_i\otimes _{{\bf Q}}R_1$ given as the tensor product of augmented rings,
together with a module, $M$,
over $R_i$.  $M$ becomes a module over $R_{i+1}$ via the following obvious
composition
$$(R_i\otimes R_1)\otimes M \fract{id\otimes id}{\ra{1.5}}
(R_i\otimes R_1) \otimes (M\otimes_{{\bf Q}}{\bf Q}) 
\fract{1\otimes T\otimes 1}{\ra{1.5}} (R_i\otimes M)\otimes
(R_1\otimes {\bf Q}) 
\fract{\mu_i\otimes \epsilon}{\ra{1.5}} M\otimes_{{\bf Q}}{\bf Q} = M.$$
Then, by change of rings,
$$Ext_{R_{i+1}}(M,{\bf Q}) ~=~ Ext_{R_i}(M,{\bf Q})\otimes
Ext_{R_1}({\bf Q}, {\bf Q}). 
\leqno{14.2}$$
Second, we will need to consider modules of the following kind:  given a module
$M$ over the ring $R$ and an exact sequence of $R$-modules
$$0 \ra{1.5} K \hookrightarrow (R)^N \ra{1.5} M \ra{1.5} 0$$
then $K$ is written $\Omega M$ and is unique up to direct sum with a
free $R$-module. 
The following result is well known.

\noindent{\bf Lemma} 14.3:  \tensl $Ext_R^{i,j}(\Omega M, N) \cong
Ext_R^{i+1, j}(M, N)$ 
for all $i\ge 1$, and
$$Ext_R^{0,*}(\Omega M, N) \ra{1} Ext_R^{1,*}(M, N)$$
is surjective for all $R$-modules $N$.\tenrm

We say $\Omega M$ is {\tensl minimal} if the natural map
$Ext_R^{0,*}(\Omega M, {\bf Q}) \ra{1} 
Ext_R^{1,*}(M, {\bf Q})$ is an isomorphism as well.  Of course, $\Omega M$
is minimal if and 
only if  the map $R^N \ra{1} M$ induces isomorphisms of $Ext^{0,*}(~,
{\bf Q})$-groups. 

Consider the following module:
$$M_1 = I \subset {\bf Q}[h_1, h_2]\leqno{14.4}$$
is the augmentation ideal, $I = (h_1,h_2)$.  Note that $M_1 = \Omega {\bf Q}$ and is
even minimal for $R = {\bf Q}[h_1,h_2]$, so
$$Ext_{{\bf Q}[h_1,h_2]}^{i,j}(M_1, {\bf Q}) = Ext_{{\bf
Q}[h_1,h_2]}^{i+1,j}({\bf Q},{\bf Q})$$ 
for all $j$ and $i\ge 0$.  Note that $Ext_{{\bf Q}[h_1,h_2]}({\bf
Q},{\bf Q}) = \ext(|h_1|^*, 
|h_2|^*)$ so 
$$Ext_{{\bf Q}[h_1,h_2]}^{i,j}(M_1, {\bf Q})~=~\cases{{\bf
Q}\{|h_1|^*, |h_2|^*\}& if 
$(i,j) ~=~ (0,2)$,\cr
{\bf Q}\{|h_1|^*|h_2|^*\}& if $(i,j) ~=~ (1,4)$.\cr}$$
It follows that a minimal resolution for $M_1$ over ${\bf Q}[h_1, h_2]$ has the
form 
$${\bf Q}[h_1,h_2]h_{1,2} \ra{1.5} {\bf Q}[h_1,h_2]\{ [h_1], [h_2]\} \ra{1.5}
M_1 \ra{1.5} 0$$
where $[h_i] \mapsto h_i$ and $h_{1,2} \mapsto h_1[h_2] - h_2[h_1]$.

Next we consider $M_1$ as a module over ${\bf Q}[h_1,h_2,h_3,h_4] =
{\bf Q}[h_1,h_2]\otimes 
{\bf Q}[h_3,h_4]$ as above and define $M_2$ as the kernel of the
surjective homomorphism 
$$\pi_2\colon {\bf Q}[h_1,h_2,h_3,h_4]\{ b_1, b_2\} \ra{1.5} M_1 \ra{1.5} 0,$$
$\pi_2(b_i) = h_i$.  Clearly $M_2 = \Omega M_1$ as a module over ${\bf
Q}[h_1, h_2,h_3,h_4]$ 
and is minimal as well.  Consequently we have
$$\eqalign{Ext_{{\bf Q}[h_1,\dots,h_4]}^{0,*}(M_2, {\bf Q})~
& =~ Ext_{{\bf Q}[h_1,\dots, h_4]}^{1,*}(M_1,{\bf Q})\cr
& =~ Ext^{2,*}_{{\bf Q}[h_1,h_2]}({\bf Q},{\bf Q})\oplus Ext_{{\bf
Q}[h_1,h_2]}^{1,*}({\bf Q},{\bf Q}) 
\otimes Ext^{1,*}_{{\bf Q}[h_3,h_4]}({\bf Q},{\bf Q})\cr
&=~ {\bf Q}^5\cr}$$
with generators $|h_1|^*|h_2|^*$, $|h_1|^*|h_3|^*$, $|h_1|^*|h_4|^*$,
$|h_2|^*|h_3|^*$, $|h_2|^*|h_4|^*$.  Likewise,
$$\eqalign{Ext^1(M_2,{\bf Q}) ~&=~Ext^2(M_1,{\bf Q})\cr
& \cong~ Ext^2_{{\bf Q}[h_1,h_2]}({\bf Q},{\bf Q})\otimes
Ext^1_{{\bf Q}[h_3,h_4]}({\bf Q},{\bf Q})\cr
&\hbox to .4in{\hfill} \oplus
Ext^1_{{\bf Q}[h_1,h_2]}({\bf Q},{\bf Q})\otimes Ext^1_{{\bf
Q}[h_3,h_4]}({\bf Q},{\bf Q})\cr 
&=~ {\bf Q}^4.\cr}$$
\indent To obtain $M_3$ regard $M_2$ as a module over ${\bf
Q}[h_1,\dots, h_6]= R_3$ 
and define $M_3 = \Omega M_2$, the kernel in the exact sequence
$$0 \ra{1.5} M_3 \hookrightarrow R^5_3 \ra{1.5} M_2 \ra{1.5} 0$$
where the map to $M_2$ takes the basis for $R^5_3$ to the five
generators of $M_2$ 
corresponding to the description above.  $M_3$ is again minimal.

In general $M_n$ is given as a module over $R_n = {\bf Q}[h_1,\dots, h_{2n}]$,
indeed is a submodule of a free $R_n$-module, and
$M_{n+1}$ is $\Omega M_n$ where $M_n$ is now regarded as a module over
$R_{n+1} = R_n\otimes {\bf Q}[h_{2n+1}, h_{2n+2}]$.

Note that for $M_1$ we have $Ext^0(M_1, {\bf Q}) = Ext^{0,2}(M_1, {\bf
Q})$ and for 
$M_2$, $Ext^0(M_2, {\bf Q}) = Ext^{0,4}(M_2, {\bf Q})$.  In other
words, for these 
modules $Ext^0$ is concentrated in a single bidegree.  Also, note that
$Ext_{R_n}( 
{\bf Q},{\bf Q}) = \ext(|h_1|^*, \dots, |h_{2n}|^*)$ has each
generator $|h_i|^*$ occurring 
in bidegree $(1, 2)$ so all the elements in $Ext^i$ occur in a single
bidegree $(i,2i)$ 
here as well.  In fact we have by a direct induction

\noindent{\bf Lemma} 14.5:  \tensl The action map
$$Ext_{R_n}^{*,*}({\bf Q},{\bf Q})\otimes
Ext_{R_n}^{0,*}(M_n, {\bf Q}) \ra{1.5} Ext^{i,j}_{R_n}(M_n,{\bf Q})$$
is surjective in all degrees and $Ext^{0,j}_{R_n}(M_n,{\bf Q}) = 0$
unless $j = 2n$. Consequently,
$$Ext^{i,j}_{R_n}(M_n,{\bf Q}) = 0$$
unless $j = 2n + 2i$.
\tenrm

As an immediate corollary we see that the $M_i$, each being minimal,
are all unique, 
since, in as much as all the generators of $M_i$ occur in the same
degree the map 
of the minimal free module onto $M_i$ is well defined up to an
isomorphism, hence 
the same is true of the kernel, $M_{i+1}$.

With this backround discussion complete we are able to state our main result.

\noindent{\bf Theorem} 14.6:  \tensl Suppose $g\ge 2$.  Then
$$Ext_{\ext_g}^{i,j}(V_g, {\bf Q}) =
\cases{{\bf Q}& if $(i,j) = (0,0)$,\cr
{\bf Q}& if $(i,j) = (1,2)$, \cr
(M_g)_{i-2, j + i-4}& if $i\ge 2$.\cr}$$
In particular, as a module over $Ext_{\ext_g}({\bf Q},{\bf Q})$
$$Ext_{\ext_g}(V_g, {\bf Q}) = {\bf Q} \oplus {\bf Q} \oplus M_g$$
where the first generator gives $Ext^{0,0}$, the second gives $Ext^{1,2}$ and
the rest is shifted by $(2, -i +4)$.
\tenrm

\noindent{\bf Example} 14.7:  We have
$$Ext_{\ext_2}(V_2, {\bf Q}) = {\bf Q} e_{0,0} \oplus {\bf Q} P_{1,2} \oplus
{\bf Q}[h_1, \dots, h_4]\{ b_1, b_2,\dots, b_5\}/{\cal R}$$
where $b_i \in Ext^{2,4}_{\ext_2}(V_2, {\bf Q})$ and ${\cal R}$ is the set
of four relations
$${\cal R} = \left\{h_5b_1,~ h_6b_1,~ h_5b_3 - h_6 b_2,~ h_5b_5 -
h_6b_4\right\}.$$ 

\noindent{\sc Proof}:  
A resolution of $V_g$ over $\ext_g$ starts in the following way
$$ 0\ra{1.5} K_g \hookrightarrow \ext_g P_{1,2} \fract{*f_g}{\ra{1.5}}
\ext_g e_{0,0} 
\ra{1.5} V_g \ra{1.5} 0$$
and we have already seen that $(K_g)_i = 0$ for $i\le g+1$.  Thus it
follows that 
$$Ext_{\ext_g}^{r,s}(V_g, {\bf Q}) = 0~\hbox{whenever $r\ge 2$ and
$s\le g+1$.}$$ 
On the other hand, since $K_g$ has $N(g,g)$ generators in $g+2$, it
follows that 
$$Ext_{\ext_g}^{2,g+2}(V_g, {\bf Q}) ~=~ {\bf Q}^{N(g,g)}.$$
{\tensl We now give an inductive analysis of the $V_g$.}  In
particular we assume 
the theorem is true for $V_{g-1}$ and begin by establishing a relation between
$V_g$ and $V_{g-1}$.  Actually, we start with $g=1$ where the theorem is not
quite true as stated, but $Ext_{\ext_1}(V_1,{\bf Q})$ is given as an extension
$$0 \ra{1.5} M_1 \ra{1.5} Ext_{\ext_1}(V_1,{\bf Q}) \ra{1.5} {\bf Q}
e_{0,0} \oplus  
{\bf Q} P_{1,2} \ra{1.5} 0$$
which turns out to be sufficient to start the induction.  There is a surjection
$$\pi_g\colon V_g \ra{1.5} V_{g-1}$$
induced by the projection $p_g\colon \ext_g \ra{1} \ext_{g-1}$,
$$p_g(e_i)~=~\cases{ e_i& for $i\le i \le 2g-2$,\cr
 0& for $i = 2g-1$ or $i = 2g$,\cr}$$
which takes $f_g$ to $f_{g-1}$.  Let $N_g$ be the kernel of $\pi_g$
and ${\cal K}_g$ 
be the kernel of $p_g$.  Then we can write
$${\cal K}_g \cong \ext_{g-1}\{ e_{2g-1}, ~e_{2g},~ e_{2g-1}e_{2g}\}$$
and
$${\cal K}_g/f_g{\cal K}_g~=~ V_{g-1}\otimes
{\bf Q}\{e_{2g-1}, ~e_{2g},~e_{2g-1}e_{2g}\}$$
as a module over $\ext_g$.  However, the surjection
$${\bar p}_g\colon {\cal K}_g/f_g{\cal K}_g \ra{1.5} N_g \ra{1.5} 0$$
is {\tensl not an isomorphism}.  In fact we have

\noindent{\bf Lemma} 14.8:  \tensl The kernel of ${\bar p}_g$ is a direct sum
of trivial modules over $\ext_g$ concentrated in degree $g+1$,
$$Ker({\bar p}_g) = Ker({\bar p}_g)_{g+1} \cong {\bf Q}^{N(g-1,g-1)}$$
as a module over $\ext_g$.
\tenrm

\noindent{\sc Proof}:  For $2\le i \le g-1$ we have
$$\eqalign{Dim(N_g)_i~&=~ \left[{2g\choose i} - {2g\choose i-2}\right]
- \left[{2g-2\choose i} - {2g-2\choose i-2}\right]\cr
&=~\left\{2{2g-2\choose i-1} + {2g-2\choose i-2}\right\} -
\left\{2{2g-2\choose i-3} 
+ {2g-2\choose i-4}\right\}\cr}$$
using the expansion
$${r\choose i}~=~ {r-2\choose i} + 2{r-2\choose i-1} + {r-2\choose i-2}$$
valid for $r \ge i \ge 2$.  This shows that ${\bar p}_g$ is an
isomorphism in this range.  Likewise,
$$\eqalign{Dim(N_g)_g~&=~ {2g\choose g} - {2g\choose g-2}\cr
&=~{2g-2\choose g} + 2{2g-2\choose g-1} - 2{2g-2\choose g-3} -
{2g-2\choose g-4}\cr 
&=~\left[ {2g-2\choose g-2} - {2g-2\choose g-4}\right] - 2\left[
{2g-2\choose g-1} - {2g-2\choose g-3}\right],\cr}$$
using the equality ${2g-2\choose g} = {2g-2\choose g-2}$, which gives
the isomorphism in dimension $g$.  However, $(N_g)_{g+1} = 0$ while
$$Dim(({\cal K}_g/f_g{\cal K}_g)_{g+1}) = Dim((V_{g-1})_{g-1})= N(g-1,g-1)$$
and the lemma follows since $({\cal K}_g/f_g{\cal K}_g)_t = 0$ for $t
\ge g+2$.\hfill\za 

Consequently we have a short exact sequence of $\ext_g$-modules
$$0 \ra{1.5} {\bf Q}^{N(g-1,g-1)} \ra{1.5} {\cal K}_g/f_g{\cal K}_g
\ra{1.5} N_g \ra{1.5} 0$$
and a long exact sequence of $Ext$-modules
$$\cdots \fract{\delta}{\ra{1.5}} Ext_{\ext_g}(N_g, {\bf Q}) \ra{1.5}
Ext_{\ext_g}({\cal K}_g, 
{\bf Q}) \fract{i^*}{\ra{1.5}}
\coprod_{i=1}^{N(g-1,g-1)}Ext_{\ext_g}({\bf Q},{\bf Q})\chi_i 
 \fract{\delta}{\ra{1.5}}\cdots\leqno{(*)}$$
where the generator, $\chi_i$ for each of the $Ext_{\ext_g}({\bf
Q},{\bf Q})$ terms occurs 
in bidegree $(0, g+1)$.

By change of rings we have
$$\eqalign{Ext_{\ext_g}({\cal K}_g, {\bf Q}) ~&\cong~
Ext_{\ext_{g-1}}(V_{g-1},{\bf Q}) 
\oplus Ext_{\ext(e_{2g-1},e_{2g})}(\{ e_{2g-1}, e_{2g},
e_{2g-1}e_{2g}\}, {\bf Q})\cr 
&=~Ext_{\ext_{g-1}}(V_{g-1},{\bf Q})\otimes {\bf
Q}[h_{2g-1},h_{2g}](b_1, b_2)/{\cal R}\cr}$$ 
where ${\cal R}$ is the relation $h_{2g}b_1 - h_{2g-1}b_2 = 0$.  By
our inductive 
assumption this implies that the generators of $Ext_{\ext_g}({\cal
K}_g,{\bf Q})$ as a 
module over $Ext_{\ext_g}({\bf Q},{\bf Q})$ occur in
bidegrees $(0,1)$, $(1,3)$, and $(2, g+2)$.  But in these bidegrees
$Ext_{\ext_g}({\bf Q},{\bf Q}) 
e_{0,g+1}$ is identically zero.  It follows that the map $i^*$ in
$(*)$ is zero so 
$$\delta\colon \coprod_{i=1}^{N(g-1,g-1)}Ext_{\ext_g}({\bf Q},{\bf
Q})\chi_i \ra{1.5} 
Ext_{\ext_g}(N_g,{\bf Q})$$
is an injection so
$$Ext_{\ext_g}(N_g,{\bf Q}) = 
\coprod_{i=1}^{N(g-1,g-1)} Ext_{\ext_g}({\bf Q},{\bf Q})e^i_{1, g+1}
\oplus Ext_{\ext_g}({\cal K}_g, {\bf Q})$$
as modules over $Ext_{\ext_g}({\bf Q},{\bf Q})$ up to a possible
extension problem. 
 
The remainder of the proof is direct.  The exact sequence
$$0 \ra{1.5} N_g \ra{1.5} V_g \fract{\pi_g}{\ra{1.5}} V_{g-1} \ra{1.5} 0$$
gives us the long exact sequence of $Ext$-groups
$$\cdots Ext^i_{\ext_g}(N_g,{\bf Q}) \fract{\delta}{\ra{1.5}}
Ext^{i+1}_{\ext_{g-1}}(V_{g-1},{\bf Q}) 
\otimes {\bf Q}[h_{2g-1}, h_{2g}] \ra{1.5} Ext_{\ext_g}(V_g,{\bf Q})
\ra{1.5} \cdots.$$ 
The map $\delta$ on the piece $Ext_{\ext_g}({\cal K}_g,{\bf Q})$ in
$Ext_{\ext_g}(N_g,{\bf Q})$ 
injects to $Ext_{\ext_{g-1}}(V_{g-1},{\bf Q})\otimes I$ where
$I\subset {\bf Q}[h_{2g-1},h_{2g}]$ 
is the augmentation ideal $(h_{2h-1}, h_{2g})$.  The resulting quotient is
$Ext_{\ext_{g-1}}(V_{g-1},{\bf Q})$ and so the calculation reduces to
the determination of the map
$$\delta\colon \coprod_{i=1}^{N(g-1,g-1)} Ext_{\ext_g}({\bf Q},{\bf
Q})\chi_i \ra{1.5} M_{g-1}$$ 
from the inductive assumption.  However, we know
$Ext_{\ext_g}^{*,j}(V_g,{\bf Q}) = 0$ 
for $*\ge 2$ and $j \le g+1$.  We also know that
$Ext_{\ext_{g-1}}^{2,g+1}(V_{g-1},{\bf Q}) = {\bf Q}^{N(g-1,g-1)}$, so
it follows that 
$\delta$ must give an isomorphism from $(1, g+1)$ to $(2, g+1)$.
Consequently, since 
$M_{g-1}$ is generated over $Ext_{\ext_{g-1}}({\bf Q},{\bf Q})$ by the
elements in this 
dimension it follows that $\delta$ is surjective and the kernel is minimal, hence
$M_g$, with the generators all occuring in bidegree $(2, g+2)$.  The
induction is complete.\hfill\za


\noindent{\bf\twelverm\S15.
The Rational Homology of $\hol{k}^*(M_g, {\bf P}^n)$ for $M_g$
Hyperelliptic}

We now use the results of \S10-\S14~ to prove the following theorem

\noindent{\bf Theorem} 15.1:  \tensl The map $H_*(\hol{k}^*(M_g,{\bf
P}^n); {\bf Q}) \ra{1} 
H_*(\map{k}^*(M_g, {\bf P}^n);{\bf Q})$ is injective for $k \ge 2g$
and $n>2$.\tenrm 

\noindent{\sc Proof}:  Consider the subquotient of the $E^1$ term
in the spectral sequence of (5.9) converging to
$H_*(QE_{v,\dots, v};{\bf Q})$ defined as the direct sum
$$\sum_{v=2s}^{g+s-1}H_*(W_{v-2s}, W_{v-2s-1})\otimes {\tilde
H}_*(S^{2s(n+1)})\otimes 
H_*(\sp{k -v}(\Sigma M_g), \sp{k-v-1}(\Sigma M_g);{\bf Q})$$
for a given $s = 0, 1, 2,\dots$
with the term
$$H_*(J(M_g), W_{g-s+1})\otimes {\tilde H}_*(S^{2s(n+1)})\otimes
H_*(\sp{k-g-s}(\Sigma M_g),\sp{k-v-s-1}(\Sigma M_g);{\bf Q})$$
added but playing a special role.  We have determined that
$$H_*(W_j,W_{j-1}) = \sum_{l=0}^{[j/2]} (V_g)_{j-2l}[M]^l$$
where $[M] \sim \sum_1^g e_{2i-1}e_{2i}$ in $H_*(J(M_g))$.  Also,
$$H_*(\sp{k}(\Sigma M_g), \sp{k-1}(\Sigma M_g)) = {\bf Q}[|e_1|,
\dots, |e_{2g}|]_k 
\oplus {\bf Q}[|e_1|,\dots, |e_{2g}|]_{k-1}|[M_g]|$$
and the (sub-quotient) $d_1$ differentials are given by $d_1(|[M_g]|) = [M_g]$,
$d_1(|e_i|) = e_i$ (including the map from the term $H_*(W_{g-s+1}, W_{g-s})
\otimes \cdots$ to the extreme term
$$H_*(J(M_g), W_{g-s})\otimes \cdots.$$
If we first apply the differential to $|[M_g]|$ we obtain the complex
$$\left\{
\sum_{v=2s}^{g+s-1}\left(V_g\right)_{v-2s}
\otimes {\bf Q}[|e_1|, \dots,
|e_{2g}|]_{k-v}\right\} \otimes {\tilde H}_*(S^{2s(n+1)}),$$
with differentials as specified as well as an appropriate quotient of
the extreme 
term with differentials mapping to it from the $v=g+s-1$-term in the sum above.

Note that the complex in brackets above is a piece of a direct summand,
$$\sum_{w=0}^{g-2}\left(V_g\right)_w\otimes
{\bf Q}[|e_1|, \dots, |e_{2g}|]_{k-2s-w}$$
of the
complex $(V_g)\otimes {\bf Q}[|e_1|, \dots, |e_{2g}|]$ which has, as its
homology the groups
$$Tor^{\ext_g}_*({\bf Q}, V_g).$$
Here, the extreme term contains the next term in the resolution as a direct
summand, but contains $H_*(J(M_g), W_{g-s-1})\otimes \cdots$ as well as
the complementary summand.
On the other hand we have already determined these $Tor$-groups in \S14~
and have
shown, in particular that $Tor^{\ext_g}_{l,m}({\bf Q},V_g) = 0$ unless
$m= m(g,l)$ is
$$m(g,l) ~=~\cases{l+g& if $l \ge 2$,\cr
2& if $l=1$,\cr
0& if $l=0$.\cr}$$
But, in the main part of the complex the $m$'s which appear are always
less than 
$m(g,l)$ so they contribute nothing by $E^2$ while at $E^2$ the term
$H_*(J(M_g), W_{g-s})\otimes \cdots$ which is a surjective image of
the stable term ${\tilde H}(J(M_g))\otimes \cdots$
contributes a quotient, consequently still a surjective image of the stable
term, and the result follows for $n >2$ by the collapsing of the stable
spectral sequence at $E^2$.\hfill\za


\centerline{\bf\twelverm References}
\vskip 10pt

\item{[ACGH]}  E. Arbarello, M. Cornalba, P.A. Griffiths, J. Harris,
{\tensl Geometry of Algebraic Curves Volume I}, Grundlehren {\bf 267},
Springer-Verlag, 1985.
\item{[BCM]}  C.-F. B\"odigheimer, F.R. Cohen, R.J. Milgram, Truncated
symmetric products and configuration spaces, {\tensl Math. Zeit.},
{\bf 214} (1993), 179-216.
\item{[BHMM]} C. Boyer, J. Hurtubise, B.M. Mann, R.J. Milgram,
The topology of the space of rational maps into generalized flag
manifolds, {\tensl Acta Math.}, {\bf 173} (1994), 61-101.
\item{[Car]} H. Cartan, {\tensl Seminaire Henri Cartan 1954-55},
exposes 2-11.
\item{[C$^2$M$^2$]} F.R. Cohen, R.L. Cohen, B.M. Mann, R.J. Milgram,
The topology of rational functions and divisors of surfaces,
{\tensl Acta. Math.}, {\bf 166} (1991), 163-221.
\item{[D]} A. Dold, Homology of symmetric products and other functors
of complexes, {\tensl Ann. of Math.}, {\bf 68} (1958), 54-80.
\item{[DT]} A. Dold, R. Thom, Quasifaserungen und unendliche
symmetrische Produkte, {\tensl Ann. of Math.}, {\bf 67} (1958),
239-281.
\item{[GH]} P. Griffiths, J. Harris, {\tensl Principles of Algebraic
Geometry}, Wiley-Interscience, 1978.
\item{[Gu]} M.A. Guest, The topology of the space of rational
curves on a toric variety, {\tensl Acta Math.}, {\bf 174} (1995),
119-145.
\item{[G]} R.C. Gunning, {\tensl Lectures on Riemann Surfaces: 
Jacobi Varieties}, Princeton U. Press, 1972.
\item{[H]} J. Hurtubise, Holomorphic maps of a Riemann surface into
a flag manifold, {\tensl J. Diff. Geom.} 1, {\bf 43} (1996), 99-118.
\item{[K1]} S. Kallel, Divisor spaces on punctured Riemann
surfaces, to appear in {\tensl transactions}.
\item{[K2]} S. Kallel, The topology of maps from a
Riemann surface into complex projective space, {\tensl preprint}.
\item{[Ku]} T. Kudo, A transgression theorem, {\tensl Mem. Fac. Sci.,
Kyusyu Univ. Ser. A}, {\bf 9} (1956), 79-81.
\item{[Mc]} I.G. Macdonald, Symmetric products of an algebraic
curve, {\tensl Topology} {\bf 1} (1962), 319-343.
\item{[Mt]} A. Mattuck, Picard bundles, {\tensl Illinois J. Math.}z, {\bf
5}(1961), 550--564.
\item{[MM1]} B.M. Mann, R.J. Milgram, Some spaces of holomorphic
maps to flag manifolds, {\tensl J. Diff. Geom.}, {\bf 33} (1991), 301-324.
\item{[MM2]} B.M. Mann, R.J. Milgram, On the moduli space of $SU(n)$
monopoles and holomorphic maps to flag manifolds, {\tensl J. Diff. Geom.}, 
{\bf 38} (1993), 39-103.
\item{[M1]} R.J. Milgram, The homology of symmetric products, {\tensl
Trans. Am. Math. Soc.}, {\bf 138} (1969), 251-265.
\item{[M2]} R.J. Milgram, The geometry of Toeplitz matrices, {\tensl
Topology}, (to appear).
\item{[S]} G. Segal, The topology of spaces of rational functions, 
{\tensl Acta. Math.}, {\bf 143} (1979), 39-72.

\vskip 40pt

\noindent Sadok Kallel\hfil\break
\noindent Centre de Recherches Math\'ematiques \hfil\break
\noindent Universit\'e de Montr\'eal \hfil\break
\noindent C.P. 6128, Succursale centre-ville \hfil\break
\noindent Montr\'eal, Quebec \hfil\break
\hbox to .3in{\hfill}Email: {\bf kallels@hans.crm.umontreal.ca}

\vskip 20pt
\noindent R. James Milgram\hfil\break
Department of Mathematics, Stanford University \hfil\break
Stanford, CA 94305, U.S.A \hfil\break
\hbox to .3in{\hfill}Email: {\bf milgram@gauss.stanford.edu}

\vfill
\eject
\bye